
\magnification=\magstep1
\baselineskip=14pt
\def \Box {\hbox{}\nobreak \hfill\vrule width 1.6mm height 1.6mm depth 0mm
\par \goodbreak \smallskip}

\rightline{\sevenrm January 11, 1994} \vskip .4cm

\centerline{\bf   BINOMIAL IDEALS }

\vskip .3cm

\centerline{\bf David Eisenbud\footnote{${}^1$}
{\rm Supported in part by the NSF.}}
\centerline{\sevenrm Brandeis
University, Waltham MA 02254}
\centerline{\sevenrm
eisenbud@math.brandeis.edu} \vskip .1cm
\centerline{and} \vskip .1cm
\centerline{\bf Bernd Sturmfels\footnote{${}^2$}
{\rm Supported in part
by the NSF and a David and Lucile Packard Fellowship.}}
\centerline{\sevenrm Cornell University, Ithaca, NY 14853}
\centerline{\sevenrm bernd@math.cornell.edu}
\vskip .4cm

\midinsert \narrower
\noindent {\bf Abstract:  }
We investigate the structure of ideals generated
by binomials (polynomials with at most two terms) and the
schemes and varieties associated to them.
The class of binomial ideals contains many classical examples
from algebraic geometry, and it has numerous applications
within and beyond pure mathematics. The ideals defining
toric varieties are precisely the binomial prime ideals.

Our main results concern primary decomposition:  If $I$ is
a binomial ideal then the radical, associated primes,
and isolated primary components of $I$ are again binomial, and
$I$ admits primary decompositions in terms of binomial primary ideals.
A geometric characterization is given for the affine algebraic
sets that can be defined by binomials. Our structural results
yield sparsity-preserving algorithms for finding
the radical and primary decomposition of a binomial ideal.

\endinsert

\beginsection Introduction.

It is notoriously difficult to deduce anything about the structure of
an ideal or scheme by directly examining its defining polynomials.
A notable exception is that of monomial ideals.
Combined with techniques for making flat degenerations
of arbitrary ideals into monomial ideals
(typically, using Gr\"obner bases),
the theory of monomial ideals becomes a useful tool
for studying general ideals. Any monomial ideal defines a scheme
whose components are coordinate planes. These objects
have provided a useful medium for exchanging information
between commutative algebra, algebraic geometry, and combinatorics.

This paper initiates the study of a larger
class of ideals whose structure can still be interpreted directly
from their generators:  binomial ideals.  By a {\it binomial}
in a polynomial ring $S=k[x_1,\dots,x_n]$ we mean a
 polynomial with at most two terms, say
$\, ax^{\alpha} + bx^{\beta}$, where $ a,b\in k $ and
$\alpha,\beta \in {\bf Z}_+^n$. We define a {\it binomial ideal}
to be an ideal of $S$ generated by binomials, and a {\it binomial
scheme} (or {\it binomial variety}, or {\it binomial algebra})
to be a scheme (or
variety or algebra) defined by a binomial ideal.
For example, it is well known that
the ideal of algebraic relations on a set of monomials
is a prime binomial ideal (Corollary 1.3).
In Corollary 2.4 we shall see that every
binomial prime ideal has essentially this form.

A first hint that there is something special about binomial
ideals is given by the following result, a weak form
of what is proved  below (see Corollary 2.4 and Theorem 6.1):

\proclaim Theorem.  The components (isolated and embedded)
of any binomial scheme in affine or projective space
over an algebraically closed field are rational varieties.

By contrast, {\sl every} scheme may be defined by trinomials,
that is, polynomials with at most three terms.  The trick is
to introduce $n-3$ new variables $z_i$ for each equation
$a_1x^{m_1}+\dots +a_nx^{m_n} = 0$ and replace this
equation by the  system of $n-2$ new equations
$$ \eqalign{
& z_1+ a_1x^{m_1}+a_2^{m_2} \quad = \quad -z_1+ z_2+a_3x^{m_3}
\quad = \quad  -z_2+ z_3+a_4x^{m_4} \quad  = \, \cdots  \cr
&  \,\,\cdots \,= \quad -z_{n-4}+ z_{n-3}+a_{n-2}x^{m_{n-2}}
\quad = \quad -z_{n-3}+a_{n-1}x^{m_{n-1}}+a_{n}x^{m_{n}}
\quad = \quad 0 .} $$

Our study of binomial ideals is partly motivated by the
frequency with which they occur in interesting contexts.
For instance, varieties of minimal
degree in projective spaces are defined by binomial equations
in a suitable system of coordinates.  More
generally, any toric variety  is defined by binomials.
(Throughout this paper we use the term ``toric variety''
to include also toric varieties that are not normal.)
Their binomial ideals are precisely the binomial prime ideals.
Sections of toric varieties by  linear subspaces defined
by coordinates or differences of
coordinates are binomial schemes.
For varieties of minimal degree such sections were studied by
Xamb\'o-Descamps [1981].

More general than coordinate rings of toric varieties
are commutative monoid algebras. An excellent general reference is the
book of Gilmer [1984], which treats these algebras
over arbitrary base rings. Gilmer shows in Theorem 7.11 that
the monoid algebras of commutative monoids are
precisely the homomorphic images of  polynomial rings
by ideals generated by
{\it pure difference binomials},
that is, polynomials $\, x^{\alpha} - x^{\beta}$, where
$\alpha,\beta \in {\bf Z}_+^n$.

Further examples generalizing toric varieties
are the {\it face rings of polyhedral complexes}
introduced by Stanley  [1987]. Geometrically, they are
obtained by gluing toric varieties along
orbits in a nice way.  They all have binomial
presentations (see Example 4.7). Some of them and their binomial sections are
geometrically interesting, for example as degenerations of special
embeddings of abelian varieties, and have played a role in the
investigations of the Horrocks-Mumford
bundle by Decker, Manolache, and Schreyer [1992].

Yet another class of algebras
with binomial defining equations is the class of {\sl Algebras
of type A} studied by Arnol\'d [1989], Korkina et al [1992] and others.
It should be possible to shed some light on their structure using
the techniques developed here.

Gr\"obner basis techniques using a total monomial order
on a polynomial ring allow the flat degeneration of an
arbitrary algebra to an algebra defined by monomial equations. Using
orders that are somewhat less strict, we sometimes get
degenerations to algebras defined
by binomial equations.  In particular, the subalgebra bases of
Robbiano and Sweedler [1990] allow one to do this in a systematic way.
The resulting degenerate
varieties may be better models of the original varieties
than those produced by a further degeneration to
varieties defined by monomials.
We hope to return to this topic in a future paper.

Complexity issues in computational algebraic geometry provide another
motivation for the study of binomial ideals. The main examples
known to attain worst case complexity for various classical problems
are binomial: these are the constructions of Mayr-Meyer [1982]
and Yap [1991] for ideal membership, Bayer-Stillman [1988] for syzygies,
Brownawell [1986] and Koll\'ar [1988] for the effective Nullstellensatz.
It has long been believed that the Mayr-Meyer schemes are so bad
because of the form of their primary decompositions.
The theory developed here provides tools for a
systematic investigation of such schemes.

Binomial prime ideals arise naturally in a variety of
settings in applied mathematics, including
dynamical systems (see e.g.~Hojevin [1992]),
integer programming (see Conti-Traverso [1991] and Thomas [1993]),
and computational statistics (see Diaconis-Sturmfels [1993]).
Within computer algebra they arise
in the extension of Gr\"obner bases theory
to canonical subalgebra bases suggested by Robbiano-Sweedler [1990],
where the role of a single S-pair is
played by an entire binomial ideal.
For real-world problems in these domains it
may be computationally prohibitive to work with
the binomial prime ideal that solves the problem exactly,
in which case one has to content oneself
with proper subideals that give approximate solutions.
Those subideals are binomial but usually not prime, so the
theory developed here may be relevant.

\vskip .1cm

We now describe the content of this paper.
To simplify the exposition, we
assume that $k$ is an algebraically closed field.
Fundamental to our treatment is the observation that
every reduced Gr\"obner basis of a binomial ideal
consists of binomials.
It follows, for example, that
the intersection of a binomial ideal and a
monomial ideal is binomial, and any projection of a binomial
scheme into a coordinate subspace has binomial closure.
Such facts are collected in Section 1, and are
used frequently in what follows.  We prove in Corollary 1.9
that the blowup algebra, symmetric
algebra, Rees algebra and associated graded algebra of a
binomial algebra with respect to a monomial ideal are
binomial algebras. This generalizes the remark
that toric blowups of toric varieties are toric.

The first step in our analysis of binomial schemes
in an affine space $k^n$ is
to decompose $k^n$ into the $2^n$ algebraic tori
interior to the coordinate planes, and  study
the intersection of the scheme with each of these.
In algebraic terms, we choose a subset
${\cal Z} \subseteq \{1,\dots,n\}$
and consider the binomial ideals in the
ring of Laurent polynomials
$$ k[{\cal Z}^{\pm}]\quad
:=\quad k[\{x_i,x_i^{-1}\}_{i\in {\cal Z}}]
\quad = \quad
S[\{x_i^{-1}\}_{i\in {\cal Z}}]/(\{x_i\}_{i\notin {\cal Z}}).
$$
These correspond to the intersections of arbitrary binomial
schemes with the tori
$$
(k^*)^{\cal Z} \quad :=\quad
\bigl\{(p_1,\dots p_n) \in k^n\ \ | \ \
p_i \neq 0{\rm\ for\ }i \in {\cal Z},\ \
p_i = 0{\rm\ for\ }i\notin {\cal Z} \bigr\}.
$$
In Section 2 we show that any binomial ideal in
$k[{\cal Z}^{\pm}]$  is a complete intersection.
In characteristic $0$ every such ``Laurent binomial ideal''
is equal to its own radical, and the algebraic set it defines
consists of several conjugate torus orbits. In
characteristic $p>0$, binomial ideals may fail to be
radical, as for example $(x^p-1)=(x-1)^p \subset k[x,x^{-1}]$,
but this failure is easy to control. We establish a
one-to-one correspondence
between Laurent binomial ideals and {\it partial characters}
on the lattice
${\bf Z}^{\cal Z}$ of monomials in $k[{\cal Z}^{\pm}]$,
where we define a partial
character $\rho$ to be a group homomorphism from
a subgroup $L_{\rho} \subseteq {\bf Z}^{\cal Z}$
to the multiplicative group $k^*$.
Properties of Laurent binomial ideals can be deduced from
arithmetic properties of the associated partial characters.
For example, the lattice $L_{\rho}$ is saturated
if and only if the corresponding Laurent binomial ideal is prime.

The next step in our theory is the study of reduced binomial
schemes.  The central result in Section 3 says that the radical of any
binomial ideal is again binomial.  We apply this in Section 4
to characterize when the intersection
of prime binomial ideals is binomial. In other words,
we determine which unions of toric varieties
are defined by binomial equations.

A serious obstacle on our road to
binomial primary decomposition lies in the fact
that if $B$ a binomial ideal and $b$ a binomial
then the ideal quotient $(B:b)$ is generally not binomial.
This problem is confronted in Section 5.
A mainspring of our theory
(Theorem 5.2) is  the description of a delicate class of
instances where these quotients are binomial.

In Section 6 we prove that the associated primes
of a binomial ideal are binomial.
Before undertaking a primary decomposition, we
pass to a ``cellular decomposition'', in which the components are
intersections of primary components having generic points
in a given cell $(k^*)^{\cal Z}$.  We then decompose the cellular
binomial ideals further: Theorem 6.4 states that the
(uniquely defined) minimal primary components are still binomial.

In Section 7 prove our main theorem:
every binomial ideal has a primary decomposition
all of whose primary components are binomial. In characteristic
$p > 0$ the result follows fairly directly from the
theory already developed, but in characteristic 0
it is much more difficult, essentially because if
$P$ is a prime binomial ideal then there will generally be
no primary binomial ideals contained in a high power of $P$.
Theorems 7.4 and 7.6 give additional information
about associated primes and primary decompositions.

In Section 8 we present some algorithms
for decomposing binomial ideals that emerge
from the general theory.  These differ markedly from the
known algorithms for primary decomposition in that they
maintain extreme sparseness of the polynomials involved.

\vskip .1cm

Having learned that the operations of primary decomposition,
radicals, projections, etc. described above take binomial ideals
to binomial ideals, the reader may think that binomiality
is preserved by many common ideal-theoretic constructions.
This is not the case; in fact, the set of ``binomial-friendly'' operations
is quite limited.  This is what makes the  main results of this paper
difficult. Here are some cautionary examples:

If $B$ is a binomial ideal and $m$ is a monomial, then
the ideal quotient $(B:m)$ is binomial (Corollary 1.7).
However, the monomial  $m$ cannot be replaced by a monomial ideal.
Even an ideal $(B:(x_i,x_j))$ need not be
binomial (Examples 1.8 and 4.6). Similarly, ideals
 $(B:b)$ for a binomial ideal $B$ and a binomial $b$
need not be binomial (Example 5.1).

Another difficulty is that very few intersections
of binomial ideals are binomial.  For example, a radical binomial
ideal can have several components, each of which must be binomial,
as stated above, but such that only certain subsets
intersect in binomial ideals.  The simplest case, in one variable,
is given by the ideal
$$
(x^d-1) \quad = \, \bigcap_{\zeta\in k,\ \zeta^d=1} (x-\zeta^m).
$$
Here the intersections of components that are again
binomial are precisely the ideals
$$
(x^{d/e}-1) \quad  = \, \bigcap_{\zeta\in k,\ \zeta^e=1} (x-\zeta^m)
$$
where $e$ divides $d$.  Our characterization of
binomial algebraic sets gives rise to
examples (such as Example 4.6) where
the intersection of the primes of maximal
dimension containing a radical binomial ideal
need not be binomial.  Given such waywardness, it
still seems to us something of a miracle that binomial ideals
have binomial primary decompositions.

\vskip .5cm

\beginsection 1. {Gr\"obner basis arguments}

Throughout this paper $k$ denotes a field
and $S:=k[x_1,\dots, x_n]$ the polynomial ring in
$n$ variables over $k$. In this section we
present some elementary facts about binomial ideals
which are proved using Gr\"obner bases..
The facts will be used frequently later on.
For Gr\"obner basics the reader may consult Buchberger [1985],
Cox, Little, and O'Shea [1992] or Eisenbud [1994].
Recall that a {\it term} is by definition a scalar
times a monomial $x_1^{i_1} x_2^{i_2} \cdots x_n^{i_n} $.

\proclaim {Proposition 1.1}.  Let $<$ be a monomial
order on $S$, and let $I \subset S$ be a
binomial ideal.
\item{(a)} The reduced  Gr\"obner basis ${\cal G}$ of $I$
with respect to $<$ consists of binomials.
\item{(b)} The normal form with respect to $<$
of any term modulo ${\cal G}$ is again a term.

\noindent {\sl Proof: } (a) If we start with
a binomial generating set for $I$,
then the new Gr\"obner basis elements
produced by a step in the Buchberger algorithm
are binomials.  \

\noindent
(b) Each step of the division algorithm modulo a set of
binomials takes a term to another term.
\Box
\vskip .1cm

One immediate application is a test for binomiality.
(Note that we are working with a fixed coordinate system.
We do not know how to test efficiently whether an ideal
can be made binomial by a linear change of coordinates.)

\proclaim {Corollary 1.2}. Let $<$ be a monomial
order on $S$.
An ideal $I \subset S$ is binomial
if and only if some (equivalently, every)
reduced Gr\"obner basis for $I$ consists of binomials.
In particular an ideal $I \subset S$ is binomial
if and only if, for every field extension $k'$ of $k$,
the ideal $k' I $ in  $k'[x_1,\ldots,x_n]$ is binomial.

\noindent {\sl Proof: } This follows from Proposition 1.1 (a)
and the uniqueness of the reduced Gr\"obner basis
with respect to a fixed monomial order ``$<$''. \Box

\vskip .1cm

Corollary 1.2 is very useful for experimentation,
since many current computer algebra systems (Axiom, Cocoa, Macaulay,
Macsyma, Maple, Mathematica,  Reduce, ...) have facilities for computing
reduced Gr\"obner bases.
The following consequence of Proposition 1.1
shows that coordinate projections
of binomial schemes are binomial:

\proclaim {Corollary 1.3}.  If $\,I \subseteq k[x_1,\dots, x_n]\,$
is a binomial ideal, then
the elimination ideal $\,I\cap k[x_1,\dots, x_r]\,$ is  a
binomial ideal for every $r \leq n$.

\noindent {\sl Proof: } The intersection is
generated by a subset of the reduced Gr\"obner basis of
$I$ with respect to the lexicographic order.
\Box

 \vskip .1cm

The projective closure is also well behaved:

\proclaim {Corollary 1.4}. If $X$ is an affine scheme in $k^n$
defined by an ideal $I$ in $S$, then the ideal in $S[x_0]$
defining the projective closure of $X$ is binomial
if and only if $I$ is binomial.

\noindent {\sl Proof: } The ideal of the
projective closure is generated by the homogenizations of the
elements in the reduced Gr\"obner basis for $I$
with respect to the total degree order.  \Box

As we have already mentioned, an intersection of
binomial ideals
is rarely binomial.  But when all but one of the ideals is generated by
monomials, or even generated by monomials modulo
a common binomial ideal, then everything is simple:

\proclaim {Corollary 1.5}. If
$I,I',J_1,\dots, J_s $ are ideals in  $S= k[x_1,\dots, x_n]$
such that $I$ and $I'$ are generated by binomials and
$J_1,\dots J_s$ are generated by monomials, then
$$
(I+I')\,\cap\, (I+J_1)\,\cap \,(I+J_2)\,\cap\,\dots \,\cap \,(I+J_s)
$$
is generated by binomials.

\noindent {\sl Proof: } Suppose first that $s=1$.
In the larger polynomial ring
$k[x_1,\dots, x_n, t\,]$ consider the binomial ideal
$L =I+tI'+(1-t)J_1$.  The claim follows from Corollary 1.3
and the formula
$ \,(I+I')\,\cap (I+J_1)\, =\, L\cap k[x_1,\dots, x_n]$.
For the general case use induction on $s$. \Box

\vskip .1cm

A slightly more subtle argument shows that there is a good
theory of monomial ideals modulo a binomial ideal.  (See
Proposition 3.4 for a further result in this direction.)

\proclaim Corollary 1.6. Let $I$ be a binomial ideal
and let $J_1,\ldots,J_s$ be monomial ideals in $S$.
\item{(a)} The intersection $\,(I+J_1)\,\cap\,\dots\,\cap\, (I+J_s)\,$
is generated by monomials modulo $I$.
\item{(b)} Any monomial in the sum $I+J_1+\cdots + J_s$
lies in one of the ideals $I+J_j$.  In particular, if
$m,m_1,\dots ,m_s$ are monomials and
$m \in I+(m_1,\dots ,m_s)$ then $m \in I+(m_i)$ for some~$i$.

\noindent {\sl Proof: }
Choose a monomial
order on $S$, and let $\cal M$ be the set of monomials
not in $in(I)$; these are called {\sl standard monomials mod I}.
The image $\overline{\cal M}$ of $\cal M$ in $S/I$ is a vector space
basis. Let $\overline J_j$ be the image of $J_j$ in $S/I$.
By Proposition 1.1 (b), each $\overline J_j$ has a vector space basis that
is a subset of $\overline{\cal M}$. It follows that the intersection
of these bases is a basis for
$\cap_j\overline J_j$, which is thus spanned by monomials.
Similarly, the union of these bases is a basis for
$ \sum_j\overline J_j$. Using Proposition 1.1 (b) again, we see
that if $m$ is a monomial in
$\sum_j (I+J_j)$ then $\overline m \in S/I$ is represented by
a standard monomial in $ \sum_j\overline J_j$, and thus belongs
to one of the $\overline J_j$, whence $m \in I+J_j$ as required.
The last statement is a special case.
\Box

Here is a central result that serves as a bridge to connect the
theory of binomial ideals in a polynomial ring
with that of Laurent binomial
ideals developed in the next section.  If $I,J$ are
ideals in a ring $R$, then we set
$\,(I:J)\ := \ \{\,f\in R \,\,| \,\, fJ \subset I\,\}$, and
$\,(I:J^\infty) \, := \, \{ \,f\in R \,| \,
f^m J \subset I\ \hbox{for}\ j \gg 0 \,\}$. If
$g\in R$, we abbreviate $(I:(g))$ to $(I:g)$.

\proclaim {Corollary 1.7}. Let $I \subset S $ be a binomial ideal,
$m_1,\dots,m_t$ monomials, and $f_1,\dots,f_t$ polynomials
such that $\sum_i f_i m_i \in I$. Let $f_{i,j}$ denote the
terms of $f_i$. For each term $f_{i,j}$, either
$f_{i,j}m_i \in I$ or there
is a term $f_{i',j'}$,
distinct  from $f_{i,j}$, and
a scalar $a \in k$ such that
$f_{i,j}m_i+af_{i',j'}m_{i'} \in I$.
In particular:
\item {(a)} For any monomial $m$ the ideal quotients
$(I:m)$ and $(I:m^\infty)$ are binomial.
\item {(b)} The first syzygies of monomials modulo
a binomial ideal are generated by binomial syzygies.

\noindent {\sl Proof: }
Choose a monomial order $>$ on $S$.
By Proposition 1.1 (b) the normal form of
$f_{i,j}m_i$ modulo $I$ is either zero or a term $m$.
If it is zero, we have $f_{i,j}m_i \in I$.
Otherwise, $m$ must cancel against a sum of terms
in the normal forms of some $f_{i'}m_{i'}$.
By Proposition 1.1 (b), these are the
normal forms of terms $f_{i',j'}m_{i'}$.
The first statement follows.

To prove (a), suppose that $f \in (I:m)$,
that is, $fm \in I$.  By the
first part of the Corollary, with $t=1$,
we may write $f$ as a
sum of binomials in $(I:m)$. Thus $(I:m)$
is generated by binomials.  Since
$(I:m^\infty) = \cup_s (I:m^s)$, the second statement
follows from the first.
Part (b) follows similarly.
\Box

\vskip .1cm

Corollary 1.7 shows that
the quotient of a binomial ideal by a
single monomial is a binomial ideal.
However, the quotient of a binomial ideal by a
monomial ideal need not be a binomial ideal, even if the
monomial ideal is generated by two variables.

\proclaim {Example 1.8}. {\sl Quotients of
binomial ideals by monomial ideals are generally not binomial}.  \rm
Let $\,I = (ax_1-ax_3,ax_2-ax_4,bx_1-bx_4,bx_2-bx_3)
\subset k[a,b,x_1,\dots,x_4]$. This ideal
is the intersection of four binomial primes
defining linear subspaces:
$$ I \quad = \quad (a, b) \, \cap \,
(a, x_1-x_4,  x_2-x_3) \, \cap \,
(b, x_1-x_3, x_2-x_4) \, \cap \,
(x_2-x_3,  x_3-x_4 ,x_1-x_4) .$$
The equidimensional part of $I$ of codimension $3$ is
$ (I:(a,b)),
$ which is the intersection of the last three of these primes. But
the homogeneous ideal
$$ (I:(a,b)) \quad = \quad \bigl(\, x_1+x_2+x_3+x_4,\,
a(x_2-x_4),\,(x_2-x_3)(x_2-x_4),\,b(x_2-x_3) \,\bigr) $$
is not a binomial ideal. For example, it contains
$x_1+x_2+x_3+x_4$ but no other linear form.
See also Example 4.6.

\vskip .1cm

Corollaries 1.3 and 1.7
give us interesting sources of
binomial algebras. For example:

\proclaim {Corollary 1.9}.
Let $B$ be a binomial ideal and $M$ a monomial ideal in $S$.
If we set $R= S/B$ and  $I=(B+M)/B \subseteq R$, then
each of the following five algebras is binomial:
the symmetric algebras
$Sym_R I$ and $Sym_{R/I} I/I^2$,
the blowup algebra
$R[zI] \subseteq R[z]$,
the Rees algebra
$R[z^{-1},zI] \subseteq R[z^{-1},z]$,
and the associated graded algebra $gr_{I} R$.

\noindent {\sl Proof: }
Let $M = (m_1,\dots,m_t)$. By Corollary 1.7
there are binomial syzygies
$\sum_j f_{i,j} m_j \equiv 0 \ (mod\ B)$
that
generate all the syzygies of
$I$ over $R$.
The symmetric algebra $Sym_R I$ may be represented as a
polynomial algebra $R[y_1,\dots,y_t]$ modulo the relations
$\sum_j f_{i,j} y_j = 0 $.  Each generator
$\sum_i f_{i,j} y_i $ is a binomial, so we see that
the symmetric algebra is binomial.
It follows that $Sym_{R/I} I/I^2 = Sym_R(I)/ISym_R(I)$
is binomial too.

The blowup algebra
$R[zI] \subseteq R[z]$
may be represented as
$R[y_1,\dots,y_t]/J$, where $J$ is the ideal of
algebraic relations satisfied over $R$ by the
elements $m_iz \in R[z]$.  The ideal $J$ is the
intersection of $R[y_1,\dots,y_t]$ with the
ideal
$$
J' = (y_1 - m_1z,\dots,y_t - m_tz) \subseteq R[y_1,\dots,y_t,z].
$$
Since $J'$ is binomial, Corollary 1.3 shows that
$J$ is binomial.  An analogous
construction with two variables $z$ and $z'$, and
an ideal
$J' = (y_1 - m_1z,\dots,y_t - m_tz, zz'-1)$
proves the statement about the Rees algebra.

The case of the associated graded algebra follows
from the cases above, since
$gr_{I} R = R[zI]/IR[zI]
= R[z^{-1},zI]/z^{-1}R[z^{-1},zI].$

\vskip .2cm

Here is another useful fact about
monomial ideals modulo binomial ideals.
The assertion is equivalent to the
existence of the special Gr\"obner basis constructed
in the proof.

\proclaim {Proposition 1.10}. Let $B$
be a binomial ideal and $M$ a monomial ideal in $S$.
If $f\in B+M$ and $f'$ is the sum of those terms of
$f$ that are not individually contained in $B+M$, then $f'\in B$.

\noindent {\sl Proof:} We may harmlessly assume
that $f=f'$, and we must show that $f\in B$.
We shall construct a special Gr\"obner basis for $B+M$.

Choose a monomial order on $S$.
Let $G$ be a Gr\"obner basis for $B$,
and let $M'$ be a set of generators for
the ideal of all monomials contained in $B+M$.  Clearly
$G\cup M'$ generates $B+M$.  We claim
 that $G\cup M'$ is a Gr\"obner basis.  By Buchberger's criterion,
it is enough to check that all s-pairs made from $G\cup M'$
reduce to zero modulo $G\cup M'$.  Now the s-pairs
made from pairs of elements of $G$ reduce to zero since $G$ is a
Gr\"obner basis.  The s-pairs made from an element of $G$ and an
element of $M'$ yield monomials that
lie in $B+M$, and that therefore
reduce to 0 through generators of $M'$. The s-pairs
made from two elements of $M'$ yield zero to begin with.
This shows that $G\cup M'$ is a Gr\"obner basis.

The normal form modulo $G\cup M'$
of a term $t$ of $f$ is, by Proposition 1.1, a monomial $m(t)$,
and our assumption implies that $m(t)$ is nonzero.
Consider the division process that reduces
$t$ to $m(t)$ by subtracting
appropriate multiples of elements of $G\cup M'$.
At each stage the remainder is a monomial.  If
this monomial were ever divisible by an element of $M'$
then it would reduce to 0. Thus the division process
can use only elements from $G$. We conclude that $f$ reduces to
zero under division by $G$, and hence $f$ lies in $B$. \Box

\vskip .5cm

\beginsection 2. Laurent binomial ideals and binomial primes

Let $k$ be a field.
We consider the ring
$$
k[{ x}^\pm] := k[{\bf Z}^n] =
k[x_1,\ldots,x_n,x_1^{-1},\ldots,x_n^{-1}]
$$
of Laurent polynomials with coefficients in $k$.
A {\it binomial} in $k[{ x}^\pm]$
is an element with at most two terms, say
$\, a { x}^{\alpha} + b { x}^{ \beta}$,
where $ a,b\in k $ and $\alpha,\beta \in {\bf Z}^n$.  A
{\it Laurent binomial ideal}  is an ideal in $ k[{ x}^\pm] $
generated by
binomials.  Note that in $k[{ x}^\pm]$ any nonzero
binomial that is  not a unit can be written
in the form $x^m-c_m$ for some $m\in {\bf Z}^n$ and
$c_m \in k^*$.

In this section we analyze Laurent binomial
ideals and their primary decompositions.
We regard $k[{x}^\pm]$ as the coordinate ring
of the algebraic torus $(k^*)^n = Hom({\bf Z}^n,k^*)$,
the group of {\it characters} of ${\bf Z}^n$.
A {\it partial character} on ${\bf Z}^n$ is a homomorphism $\rho$
from a sublattice $L_{\rho}$ of ${\bf Z}^n$ to the
multiplicative group $k^*$. Whenever we speak of a
partial character $\rho$, we mean the pair consisting of the
map $\rho$ and its domain $L_\rho \subseteq {\bf Z}^n$.
Given a partial character $\rho$,
we define a Laurent binomial ideal
$$
{ I}(\rho) \quad :=  \quad ( {x}^{m} - \rho({m})
\,:\, {m} \in L_{\rho} ) .
$$
We shall see that all Laurent binomial ideals are of this form.

The algebraic set $Z(I(\rho))$
of points in $(k^*)^n = Hom({\bf Z}^n,k^*)$ where all
the elements of $I(\rho)$ vanish is precisely the set of characters
of ${\bf Z}^n$ that restrict to $\rho$ on $L_{\rho}$.
If $k$ is algebraically closed, then $Z(I(\rho))$
is nonempty for any partial character
$\rho$. This follows from the Nullstellensatz, or
from the fact that the group $k^*$ is divisible.

If $L$ is a sublattice of ${\bf Z}^n$,
then the  {\it saturation} of $L$ is the lattice
$$ Sat(L) \quad := \quad \{\, {m} \in {\bf Z}^n \,\,|\,\, d{m}
\in L\,\, {\rm for\, some}\,\, d \in {\bf Z} \,\}. $$
The group $Sat(L)/L$ is finite. We say that
$L$ is {\it saturated} if $L=Sat(L)$.

\proclaim Theorem 2.1.  Let
$k[{x}^\pm]$ be a Laurent polynomial ring
over a field $k$.
\item{(a)} For any proper Laurent binomial ideal
$I \subseteq \, k[{x}^\pm]$
there is a unique partial character $\rho$ on
${\bf Z}^n$  such that $I = I(\rho)$.
\item{(b)} If $m_1,\dots, m_r$ is a basis of the lattice
$L_{\rho}$, then the binomials
$$ x^{m_1}-\rho(m_1),\dots,x^{m_r}-\rho(m_r) $$
generate $I(\rho)$ and form a regular sequence in $k[x^{\pm}]$.
In particular
$$ codim(I(\rho))=rank(L_{\rho}). $$
\vskip .1cm
Now assume that $k$ is algebraically closed.
\item{(c)} The ideal $I(\rho)$ is prime if and only if
$L_{\rho}$ is saturated. In this case $Z(I(\rho))$ is the
orbit of the point  $(\tilde \rho(e_1),\dots,\tilde\rho(e_n)) $
under the group of characters of ${\bf Z}^n$ that are
trivial on $L_{\rho}$, where $\tilde \rho$ is any
extension of $\rho$ to ${\bf Z}^n$.
\item{(d)} Let $char(k) =p \geq 0$.
Suppose that $\rho$ is a partial character on ${\bf Z}^n$ and
$L_{\rho} \subseteq L \subseteq {\bf Z}^n$ are lattices
with $L/L_{\rho}$ finite of order $g$.
If $g$ is relatively prime to $p$,
then there are $g$ distinct characters
$\rho'$ on $L$ that are extensions of
$\rho$ on $L_{\rho}$,
and
$$
I(\rho) \quad  = \,\bigcap_{\rho'\,extends\, \rho\, to\,
L} I(\rho'). $$
If $g$ is a power of $p$,
then there is a unique extension $\rho'$ of $\rho$ to $L$,
and $k[{x}^\pm]/I(\rho)$ has a filtration by $k[{x}^\pm]$-modules
$$
k[{x}^\pm]/I(\rho)= M_0 \supset M_1 \supset \dots \supset M_g = 0
$$
with  successive quotients $M_i/M_{i+1} \cong k[{x}^\pm]/I(\rho')$.
\item{ }

\vskip -.3cm

\noindent {\sl Proof: }
(a)  Any proper binomial ideal $I$ in
$k[x^\pm] $ is generated by its elements of the form
$\,x^m - c$ for $m \in {\bf Z}^n$ and $c \in k^*$.
Let $L$ be the
subset of
${\bf Z}^n$
consisting of those
$m$
that appear.  Since $I$ is proper,
$c_m$ is uniquely determined by $m$.
{}From the basic formula
$$
x^{m+m'}-cd \quad = \quad (x^m-c)x^{m'}-c(x^{m'}-d)  \eqno (2.1)
$$
we see that if $\,{x}^{{m}} - c_{{m}} $
and
$\,{x}^{{m'}} - c_{{m'}} $ are in ${ I}$,
then
so is
$
\,{x}^{{m} + {m'}} - c_{{m}}c_{{m'}} \,
$
while if
$\,{x}^{{m+m'}} - c_{{m}}c_{{m'}} $
and
$\,{x}^{{m}} - c_{{m}} $ are in ${ I}$, then so is
$\,{x}^{{m'}} - c_{{m'}}$.
Hence $L$ is a
sublattice of ${\bf Z}^n$, the map $\rho\,:L \to k^*$ taking
$m$ to $c_m$ is a character, and $I = I(\rho)$.

For the uniqueness part of (a) we shall show that
if a binomial $x^u - c_u $ lies in $I(\rho)$
then $ u \in L_\rho$ and $c_u = \rho(u)$.
We write $k[x^\pm]$ as the quotient of the polynomial
ring $T := k[y_1,\ldots,y_n,z_1,\ldots,z_n]$ modulo
the binomial ideal $(y_i z_i -1 : i=1,\ldots,n)$.
If $I'(\rho)$ denotes the preimage of $I(\rho)$
in $T$, then $I'(\rho)$ is generated by
the set
$$\{ y^a z^b -
\rho(a\!-\!b\!-\!c\!+\!d) \cdot y^c z^d \,: \,a,b,c,d \in {\bf N}^n,
a-b-c+d\in L_\rho \}.\eqno (2.2)
$$
In fact, this set is a Gr\"obner basis for  $I'(\rho)$
with respect to any monomial order on $T$, by Proposition 1.1.
If $x^u - c_u$ lies in $I(\rho)$,
and we write $u_+, u_-$ for the positive and
negative parts of $u$, so that
$u= u_+-u_-$, then the normal form of
$y^{u_+} z^{u_-} $  modulo this
Gr\"obner basis is the constant $c_u$.
Each polynomial in the reduction sequence
is a term of the form
$\,\rho(a\!-\!b\!-\!c\!+\!d)\cdot y^{u_+ - a + c} z^{u_- - b + d}\,$
where $a-b-c+d\in L_\rho $.
This proves that $u \in L_\rho$ and $\rho(u) = c_u$.

(b)  Formula (2.1)
shows that any set of additive generators
$\{m_i\}$ of $L_{\rho}$ gives rise to a set of
generators $x^{m_i}-\rho(m_i)$ of the ideal $I(\rho)$.

If $m_1,\dots,m_r$ are linearly independent elements
that span $L_{\rho}$
it remains to show that
$$
x^{m_1}-\rho(m_1),\dots,x^{m_r}-\rho(m_r)
$$
is a regular sequence.  By
induction on $r$ we may suppose that the first $r-1$
binomials
form a regular sequence.  In particular all the
associated primes of
the ideal they generate have codimension $r-1$.
Thus it suffices to show that
the ideal $I(\rho)$
has codimension $r$.

Let $L$ be the saturation of $L_{\rho}$.  We may write
${\bf Z}^n = L\oplus L'$ for some lattice $L'$, so
$k[{\bf Z}^n]/I(\rho) = k[L]/I(\rho)\otimes k[L']$,
which is a Laurent polynomial ring in $n-r$ variables over
$k[L]/I(\rho)$.  Thus
it suffices to show that $k[L]/I(\rho)$ has dimension 0,
or using the Nullstellensatz, that the set of
characters of
$L$  that are extensions of $\rho$ is finite.

{}From the exact sequence
$$
0\to L_\rho \to L \to L/L_\rho \to 0
$$
we see that any
two characters of $L$ restricting to $\rho$ on $L_\rho$
differ by a character of the finite group
$L/L_{\rho}$.  Since $k$ is a field, its subgroup of
elements of any given finite order is cyclic, and in particular
finite.  Thus $Hom(L/L_{\rho},k^*)$ is finite as required.

(c) Suppose that
$L = L_{\rho}$ is saturated. Writing ${\bf Z}^n = L\oplus L'$
as before we get
$$
k[{\bf Z}^n]/I(\rho)\,\,\, = \,\,\, k[L]/I(\rho)\otimes_k
k[L'] \,\,\,= \,\,\, k\otimes_k k[L'] \,\, =\,\, k[L'].
\eqno (2.3) $$
This is a domain, hence $I(\rho)$ is prime.

Conversely, suppose $I(\rho)$ is prime.  If $m\in {\bf Z}^n$ and
$dm\in L_{\rho}$ then
$$
x^{dm}-\rho(dm)
\quad = \quad \prod_{i=1}^d (x^m-\zeta^i\rho(m))
\quad \in\quad I(\rho)
$$
where $\zeta$ is a generator of the
group of $d^{\rm th}$ roots of unity in $k$.
Thus one of the factors $x^m-\zeta^i\rho(m)$ belongs
to $I(\rho)$,  and we see that $m \in L_{\rho}$ by the
uniqueness  statement of part (a).  Thus $L_{\rho}$
is saturated.

If $L = L_\rho$ is saturated, then
the group of characters of ${\bf Z}^n = L\oplus L'$ that are trivial
on $L$ may be identified with
 the group of characters of $L'$. The last
statement of (c) now comes from the identification of
$Z(I(\rho))$ with the set of characters extending $\rho$.

\vskip .1cm

(d) Both statements reduce immediately to the case where
$L/L_{\rho}$ is a cyclic group of prime order $q$.
Diagonalizing a matrix for the inclusion $L_{\rho}\subset L$ we
may choose a basis $m_1,\dots,m_r$ of $L$ such that
$L_{\rho}$ has the basis $m_1,\dots,m_{r-1}, qm_r$.
For any extension $\rho'$ of $\rho$ to $L$, the element
$\rho'(m_r)$ is a $q^{\rm th}$ root of $\rho(qm_r)$.  If
$c \in k^*$ is one such $q^{\rm th}$ root and
we let
$$ J=(x^{m_1}-\rho(m_1),\dots,x^{m_{r-1}}-\rho(m_{r-1})) $$
then each of the  ideals $I(\rho')$ has the form
$ I(\rho')=J+(x^{m_{r}}-\zeta c) $
for some $q^{\rm th}$ root of unity $\zeta$, while
$I(\rho) = J +(x^{qm_{r}}-c^q)$.

If $q\neq p$, then there are $ q$ distinct
$q^{\rm th}$ roots of unity in $k$.
If $\zeta$ and $\zeta'$ are two of them then
$ I(\rho')=J+(x^{m_{r}}-\zeta c)$ and
$ I(\rho')=J+(x^{m_{r}}-\zeta' c)$ together generate the unit
ideal.  Thus in the ring $R := k[x^{\pm}]/J$ the intersection
of these ideals is equal to their product, and we get
$$ I(\rho)/J \,\,= \,\, (x^{qm_{r}}-c^q)R
\,\,= \,\,\prod_\zeta
(x^{m_{r}}-\zeta c)R \,\, = \,\,\bigcap_{\zeta}
(x^{m_{r}}-\zeta c)R \,\,= \,\, \bigcap_{\rho'} I(\rho')/J.
$$
It follows that
 $ I(\rho)=\bigcap_{\rho'} I(\rho')$ as required.

On the other hand, if $q = p$ then
$\zeta = 1$ and $x^{qm_{r}}-c^q=(x^{m_{r}}-c)^q$.  By part
(b), the element $x^{m_{r}}-c$ is a nonzerodivisor
modulo $J$.
Therefore in the filtration
$$
k[{x}^\pm] \supset I(\rho') = J+(x^{m_{r}}-c)
\supset J+(x^{m_{r}}-c)^2 \supset
\dots \supset J+(x^{m_{r}}-c)^p = I(\rho),
$$
the successive quotients are isomorphic to $k[{x}^\pm]/I(\rho')$.
 Reducing
modulo $I(\rho)$, we get a filtration of $k[{x}^\pm]/I(\rho)$ with the
desired properties.\Box

\vskip .2cm

Using Theorem 2.1 we can describe the primary decomposition
and radical of a Laurent binomial ideal in  terms of operations
on integer lattices.  If $L$ is a sublattice of ${\bf Z}^n$,
and $p$ is a prime number,
we define $Sat_p(L)$ and $Sat'_p(L)$ to be
the largest sublattices of $Sat(L)$ such that
$Sat_p(L)/L$ has order  a power of $p$ and
$Sat'_p(L)/L$ has order relatively prime to $p$.
(These can be computed by diagonalizing a matrix
representing the inclusion of $L$ in ${\bf Z}^n$.)  We
adopt the convention that if $p=0$ then $Sat_p(L) = L$ and
$Sat'_p(L) = Sat(L)$.

If $\rho$ is a partial character, we define the {\it saturations}
of $\rho$ to be the characters
$\rho'$ of $Sat(L_{\rho})$
that restrict to $\rho$ on $L_{\rho}$, and we say that
$\rho$ is saturated if $L_\rho$ is saturated.

\proclaim Corollary 2.2. Let $k$ be an algebraically
closed field of characteristic $p \geq 0$.
Let $\rho$ be a partial character.
Write $g$ for the order of $Sat'_p(L_{\rho})/L_{\rho}$.
There are $g$
distinct characters $\rho_1,\dots,\rho_g$ of
$Sat'_p(L_{\rho})$ extending $\rho$ and for each $j$ a unique
character $\rho_j'$ of $Sat(L_{\rho})$ extending $\rho_j$.
There is a unique partial character of $\rho'$ of
$Sat_p(L_{\rho})$ extending $\rho$.
The radical, associated primes, and minimal primary
decomposition of
$I(\rho) \subseteq k[x^\pm]$ are:
$$\eqalign{
\sqrt{I(\rho)} \quad &= \quad I(\rho')\cr
Ass(S/I(\rho)) \quad &= \quad \{I(\rho_j')\ |\
j=1,\dots,g\}\cr I(\rho) \quad &= \quad \bigcap_{j=1}^g
I(\rho_j), }$$
and $I(\rho_j)$ is $I(\rho_j')$-primary.
In particular, if
$p = char(k) = 0$ then $I(\rho)$ is a radical ideal.
The associated primes $I(\rho_j')$ of $I(\rho)$ are all minimal
and have the same codimension $rank(L_{\rho})$.
The geometric multiplicity of each
primary component $I(\rho_j)$ is the order of the group
$Sat_p(L_{\rho})/L_{\rho}$.

\noindent {\sl Proof:} For every prime $q \neq p$ and every
integer $d \geq 0$ the subgroup of $k^*$ of elements of order
$q^d$ is cyclic of order $q^d$, while the subgroup of $k^*$ of
elements of order $p^d$ is trivial.  This implies that
 there is a unique extension $\rho'$
of $\rho$ to $Sat_p(L_{\rho})$, exactly $g$
extensions $\rho_j$ of $\rho$ to $Sat'_p(L_{\rho})$,
and a unique extension $\rho_j'$
of $\rho_j$ to $Sat(L_{\rho})$.
Since  $Sat(L_{\rho})/L_{\rho}$ is finite,
the rank of $Sat(L_{\rho})$ is the same
as that of $L_{\rho}$.

By Theorem 2.1 (b) and (c),
each $I(\rho_j')$ is a prime ideal of codimension
$=rank(L_{\rho})$.  By the first part of
Theorem 2.1 (d) we have $I(\rho') = \cap_j I(\rho_j')$,
so $I(\rho')$ is a radical ideal.  The second part of
Theorem 2.1 (d) shows that $k[x^{\pm}]/I(\rho)$ has a finite
filtration whose factors are isomorphic to $k[x^{\pm}]/I(\rho')$,
so that $I(\rho')$ is nilpotent mod $I(\rho)$.  This shows that
$I(\rho')$ is the radical of $I(\rho)$.

The equality $I(\rho) = \bigcap_{j=1}^g I(\rho_j)$
follows directly from the first part of Theorem 2.1 (d).
Thus to establish the assertions about associated primes
and primary decomposition,
it suffices to show that each $I(\rho_j)$ is
$I(\rho_j')$-primary of geometric multiplicity
$card(Sat_p( L_\rho)/L_\rho)$. Applying the second part of
Theorem 2.1 (d), we see that $k[x^{\pm}]/I(\rho_j)$ has a
filtration of length $g$ whose successive quotients are
all isomorphic to $k[x^{\pm}]/I(\rho_j')$. Both the fact
that $I(\rho_j)$ is primary and the assertion about the
geometric multiplicity follow.
\Box

\vskip .1cm

The results of Theorem 2.1 can be transferred to certain
affine binomial ideals.  As in the proof of Theorem 2.1
(a), we let $m_+, m_-\in {\bf Z}_+^n$ denote the
{\it positive part} and {\it negative part} of a
vector $m\in {\bf Z}^n$. Given a partial character
$\rho$ on ${\bf Z}^n$, we define the ideal
$$ \, { I}_+(\rho) \quad :=  \quad (\{ {x}^{m_+} - \rho({m})x^{m_-}
\,:\, m \in L_{\rho}\} )\quad
\hbox{in} \quad  S = k[x_1,\ldots,x_n]. \eqno (2.4)$$

\proclaim {Corollary 2.3}. If $I$ is a binomial ideal
in $S = k[x_1,\ldots,x_n]$ not containing any
monomial, then there is a unique partial character
$\rho$ on ${\bf Z}^n$ such that
$\, (I:(x_1\cdots x_n)^{\infty})\, = \, I_+(\rho)$.
The generators of $I_+(\rho)$ given in (2.4) form
a Gr\"obner basis for any monomial order on $S$.
The binomial ideals of the form $I_+(\rho)$
are precisely those whose associated points are
off the coordinate hyperplanes.  If $k$ is algebraically
closed, then all the statements of
Corollary 2.2 continue to hold if we replace each
$I(-)$ by $I_+(-)$.

\noindent {\sl Proof: }
The ideal $\, (I:(x_1\cdots x_n)^{\infty}) \, $
is equal to  $\,I \cdot k[x^{\pm} ] \, \cap \, S $,
the contraction from the Laurent polynomial ring.
By Theorem 2.1 (a), there exists a unique partial character $\rho$
such that $\,I \cdot k[x^{\pm} ] \, = \,  I(\rho) \cdot k[x^{\pm} ]$.
The map $S\to k[x^\pm]$ may be
factored through the ring $T$ as in
the proof of Theorem 2.1 (a).  With
$I'(\rho)$ defined as in that proof, we have
$\,I \cdot k[x^{\pm} ] \, \cap \, S = I'(\rho) \, \cap
\,S$. Since the elements in the set (2.2) form a Gr\"obner
basis with respect to any monomial order on $T$,
the elements in this set not involving the variables $y_i$
form a Gr\"obner basis of
$\,I \cdot k[x^{\pm} ] \, \cap \, S $.  These are
exactly the given generators of $I_+(\rho)$.

The third statement holds because
an ideal in $S$ whose associated points are off the
coordinate hyperplanes is contracted from $k[x^{\pm}]$.
The fourth statement follows at once. \Box

\vskip .1cm

Consider a $k$-algebra homomorphism from $S =
k[x_1,\dots,x_n]$ to the Laurent polynomial ring
$k[t^\pm] := k[t_1,t_1^{-1},\dots,t_r,t_r^{-1}]$
which sends each variable $x_i$ to a
monomial $c_i t^{a_i}$. Its kernel $P$ is a prime
ideal, which is generated by binomials.
The variety defined by $P$ in $k^n$ is a
(not necessarily normal)  affine toric variety.
For details on toric varieties
and their ideals see  Fulton [1993],
Sturmfels [1991], and the references given there.
Corollary 2.3 implies that the class of
{\it toric ideals} is the same as the class
of binomial prime ideals.

\proclaim {Corollary 2.4}. Let $k$ be an algebraically
closed field, and let
$P$ be a binomial ideal in $S=k[x_1,\dots,x_n]$.
Set
$\,\{y_1,\dots,y_s\} := \{x_1,\ldots,x_n\}\,\cap \, P\,$
and let
$\,\{z_1,\dots,z_t\} := \{x_1,\ldots,x_n\}\, \setminus P$.
The ideal $P$ is prime if and only if
$$
P = (y_1,\dots,y_s)+I_+(\rho)
$$
for a saturated partial character $\rho$  in
the lattice ${\bf Z}^t$ corresponding to $z_1,\dots,z_t $.
In this case, the prime $P$ is the kernel of
a ring homomorphism
$$ k[y_1,\dots,y_s,z_1,\dots,z_t] \,\to
\, k[t^{\pm}]\,, \,\,\, y_i \,\mapsto \,0 \,,\,\,\,
z^m \mapsto \tilde\rho( m) t^{\bar m}, \eqno (2.5) $$
where $\bar m \in {\bf Z}^t/L_\rho$ denotes
the image of $m \in {\bf Z}^t$, the
group algebra of ${\bf Z}^t/L_\rho$
is identified with a Laurent polynomial ring $k[t^{\pm}]$,
and $\tilde\rho$ is any extension
of $\rho$ to ${\bf Z}^t$.

\noindent {\sl Proof:  } We must prove the
``only if''-direction. Given a binomial prime $P$,
consider the binomial prime
$P/(y_1,\ldots,y_s)$ in $k[z_1,\dots,z_t]$. Modulo
this prime each $z_j$ is a nonzerodivisor. By Corollary 2.3,
we may write $P/(y_1,\ldots,y_s) = I_+(\rho)$.
Since $Pk[z^{\pm}] = I(\rho)$ is prime,
Theorem 2.1 (c) shows that $\rho$ is saturated.
For the proof of the second statement consider the
surjective homomorphism $\,k[z^\pm]
\rightarrow k[t^{\pm}],\,
z^m \mapsto \tilde\rho( m) t^{\bar m}$.
Its kernel  obviously contains $Pk[z^{\pm}]$,
and since $Pk[z^{\pm}]$ is a prime of codimension
$rank(L_\rho) = dim(k[z^\pm]) - dim(k[t^\pm])$,
the kernel is precisely $P$.  Since $P$ is the preimage of
$Pk[z^{\pm}]$ in $S$, we conclude that $P$ is the kernel of
the  composite map $S \to k[z^\pm] \to k[u^{\pm}]$,
which coincides with (2.5).
\Box

 \vskip .5cm

\beginsection 3. The radical of a binomial ideal

The {\it radical} of an ideal $I $ in $S =
k[x_1,\ldots,x_n]$ is
$\,\sqrt{I} \,:= \,\{\,f\in S \,| \,f^d\in I \,\,\hbox{for} \, \, d
\gg 0\}$. In this section we show that the family of
binomial ideals is closed under taking radicals.

\proclaim Theorem 3.1. Let $I\subseteq S  = k[x_1,\ldots,x_n]$
be an ideal.  If $I$ is binomial then $\sqrt{I}$ is binomial.

In the special case where $I$ is generated by
pure difference binomials (monomial minus monomial), this
result was proved using different methods
by Robert Gilmer [1984, section 9];
Gilmer's results show that the radical
is again generated by pure difference binomials, and prove
a similar statement for the case of an arbitrary base ring.

Our proof works by an induction on the number of variables, and
an application of the Laurent case treated in the previous
section.  For this we use:

\proclaim Lemma 3.2. Let $R$ be any commutative ring, and
let $x_1,\ldots,x_n\in R.$ If
$I$ is any ideal in $R$, then the radical of $I$
satisfies the relation
$$ \sqrt{I} \quad = \quad
\sqrt{ \bigl( I : (x_1 \cdots x_n)^\infty \bigr)} \,
\,\cap \, \sqrt{  I + (x_1) } \,\,\cap \,\,
\cdots \,\,\cap \, \sqrt{ I + (x_n) }. \eqno
(3.1) $$

\noindent {\sl Proof: } The right hand side clearly
contains $\sqrt{I}$.
It suffices to show that every prime $P$ containing $I$
contains one of the ideals on the right hand side.
If $\,( I : (x_1 \cdots x_n)^\infty) ) \subseteq P$ we are
done. Otherwise,
$\,f \cdot (x_1 \cdots x_n)^d \in I \subset P\,$ for
some integer $d$ and some $f \in R \setminus P$.
This implies $x_i \in P$ for some $i$.
Thus $P$ contains $I + (x_i)$ as required. \ \Box

\proclaim Lemma 3.3. Let $I$ be a binomial ideal in
$S = k[x_1,\ldots,x_n]$. Set $S' = k[x_1,\ldots,x_{n-1}]$.
If $I' = I \cap S'$, then
$I + (x_n)$ is the sum of $I' S + (x_n)$ and
an ideal generated by monomials in $S'$.

\noindent {\sl Proof: }
Every binomial that involves $x_n$ is either
contained in $(x_n)$ or is congruent modulo $(x_n)$
to a monomial in $S'$. Thus all generators of $I$
which are not in $I'$ may be replaced by monomials
in $S'$ when forming a generating set for $I + (x_n)$.
\ \Box

\proclaim Proposition 3.4.
Let $I$ be a binomial ideal in $S$.
If $M$ is a monomial ideal, then
$\sqrt{I+ M} = \sqrt{I} + M_1$ for some
monomial ideal
$M_1$.

\noindent {\sl Proof: }  We may suppose that
$I= \sqrt{I}$.
We apply Lemma 3.2 to the ideal $I+M$.
If $M = (0)$ there is
nothing to prove, so we may assume that $M$
actually contains a monomial.
In this case $((I+M): (x_1\cdots x_n)^\infty) = S$,
and Lemma 3.2 yields
$\, \sqrt{I+M} \, = \, \cap_{i=1}^n \sqrt{I+M + (x_i)}$.
By Corollary 1.5, it suffices to show that
the radical of $I+M + (x_i)$ is the sum of $I$
and a monomial ideal.

For simplicity let $i=n$  and write
$S' = k[x_1,\ldots,x_{n-1}]$. Since $I$ is radical,
the ideal $I' = I \cap S$
is radical as well. By Lemma 3.3,
$I + M + (x_n) = I'S + JS + (x_n)$
where $J$ is a monomial ideal in $S'$.
By induction on $n$, the radical of $I' + J$ in $S'$
has the form $I' + J_1$, where $J_1$ is a monomial
ideal of $S'$. This implies the following identity
of ideals in $S$:
$$ \sqrt{I + M + (x_n)}
\quad = \quad \sqrt{ I'S + JS + (x_n) } \quad
=  \quad I' S + J_1 S + (x_n) \quad = \quad
I  \,+\, J_1 S + (x_n) .$$
\vskip -0.7cm
\Box

\noindent {\sl Proof of Theorem 3.1.  }
We proceed by induction on $n$, the result being
trivial for $n=0$. Let $I$ be a binomial ideal in $S$.
Let $I_j := I\cap S_j$
where $S_j = k[x_1,\ldots,x_{j-1},x_{j+1},
\ldots,x_n]$.
By induction we may  assume that the radical of
each $I_j$ is binomial.  Adding these binomial ideals to
$I$, we may assume that each $I_j$ is radical to begin with.

We shall use the formula (3.1) in Lemma 3.2 for $\sqrt{I}$.
The ideal $\,\sqrt{\bigl( I : (x_1 \cdots x_n)^\infty \bigr)}\,$
is binomial by Corollaries 1.7, 2.2 and 2.3, and we can write
it as $I+I'$ for some binomial ideal $I'$.
To show that the intersection in formula
(3.1) is binomial we use Corollary 1.5. It suffices
to express $\sqrt{I + (x_j)} $ as the sum of $I$ and a
monomial ideal.
By Lemma 3.3, we can write
$\,I + (x_j)\,= \, I_j S + J S + (x_j)$,
where $J$ is a monomial ideal in $S_j$.
The radical of  $\,I + (x_j)\,$ equals the
radical of $ \, I_j S + JS + (x_j)$.
However, $I_j $ is a radical binomial ideal
by our inductive assumption made above.
We can apply Proposition 3.4
with $M = JS + (x_i)$ to see that there exists
a monomial ideal $M_1$ in $S$ such that
$$\,\sqrt{I + (x_j)} \,\,=
\,\,\sqrt{ \,{I_j} S + J S + (x_j) \,}
\,\,= \,\, {I_j} S + M_1 \,\, = \,\,I + M_1 .
$$
\vskip -0.7cm
\ \Box

\vskip .3cm

\proclaim Example 3.5. \  (Permanental ideals) \ \ \ \rm
We do not know how to tell whether a binomial ideal is radical
just from the shape of a generating set. As an example consider
the ideal $P_{m,n}$ generated by the $2 \times 2$-subpermanents
$\,x_{ij} x_{kl} + x_{il} x_{kj}\, $ of an $m \times n$-matrix
$(x_{ij})$ of indeterminates over a field $k$ with $char(k) \not= 2 $.
If $m \leq 2 $ or  $n \leq 2 $ then $P_{m,n}$ is a radical ideal.
(This can be shown using the technique in Proposition 4.8).
For instance, we have
$$ \eqalign{ P_{2,3} \quad = \quad\,\,
& (x_{11}, x_{12},x_{13}) \,\, \cap \,\,
(x_{21}, x_{22},x_{23}) \,\, \cap \,\,
( x_{11} x_{22} + x_{12} x_{21}, x_{13}, x_{23}) \cr
& \,\,\cap \,\, ( x_{11} x_{23} + x_{13} x_{21}, x_{12}, x_{22})\, \cap \,
  (x_{12} x_{23} + x_{13} x_{22}, x_{11}, x_{21}) .\cr }  $$
However, if $m, n \geq 3$ then $P_{m,n}$ is not radical:
$\,x_{11}^2 x_{22} x_{33} \in P_{m,n} \,$
but $\,x_{11} x_{22} x_{33} \not\in P_{m,n} $.  Of course if
the plus signs in the generators of $P_{m,n}$ are changed
to minus signs we get a determinantal ideal that is
prime for every $m$ and $n$.

\vskip .6cm

\beginsection 4. Binomial algebraic sets

We next characterize intersections
of prime binomial ideals that are
generated by binomials.  The result is best stated
geometrically.  For this purpose we define
an algebraic set to be a reduced affine algebraic
scheme over $k$.  (Alternately, one may work with
ordinary algebraic sets defined by equations
with coefficients in $k$ but having points
with coordinates in some fixed algebraic closure of $k$;
or one may simply restrict to the case where
$k$ is algebraically closed.) \  By Theorem 3.1, an
algebraic set is cut out by binomials
set-theoretically if and only if its
ideal is  generated by binomials.
Such a set is called a {\it binomial algebraic set}.

We decompose affine $n$-space $k^n$
into tori corresponding to the
$ 2^n $ coordinate flats
$$ (k^*)^{\cal Z} \quad :=\quad
\bigl\{\,(p_1,\ldots, p_n) \in k^n\ \ | \ \
p_i \neq 0{\rm\ for\ }i \in {\cal Z},\ \
p_i = 0{\rm\ for\ }i\notin {\cal Z}\, \bigr\}, \eqno (4.1)
$$  where  ${\cal Z}$ runs over all subsets of
$\{1,\ldots,n\}$. We shall refer to the
tori $(k^*)^{\cal Z}$ as {\it coordinate cells}.
The closure of a coordinate cell $(k^*)^{\cal Z}$ in $k^n$
is defined by the ideal
$$ M({\cal Z}) \,\,\, := \,\,\,
(\{x_i\ |\ \ i\notin {\cal Z}\})
\quad \quad \hbox{in} \,\quad
 S \,\, = \,\,k[x_1,\dots,x_n] .$$
The coordinate ring of $(k^*)^{\cal Z}$ is the
Laurent polynomial ring
$$ k [ {\cal Z}^\pm ] \quad := \quad
k\, [\,\{x_i,x_i^{-1}\}_{ i\in {\cal Z}}\, ] . $$
There is a {\it coordinate projection}
$ (k^*)^{{\cal Z}'}\to (k^*)^{\cal Z}$
whenever ${\cal Z}\subseteq {\cal Z}'\subseteq
\{1,\ldots,n\}$. It is   defined by
setting  all those coordinates not in
${\cal Z}$ to zero.

If $X$ is any subscheme of $k^n$,
corresponding to an ideal  $I\subseteq S $
then the closure of the
intersection of $X$ with the coordinate cell
$(k^*)^{\cal Z}$ corresponds to the ideal
$$ I_{\cal Z}\quad :=\quad
\biggl( (I + M({\cal Z})) \, :\,
(\prod_{i\in {\cal Z}} \! x_i \,)^{\infty}\biggr).\eqno
(4.2) $$
This ideal can be identified with the image of $I$
in $k[{\cal Z}^\pm]$.
If $I$ is radical, then it is easy to see
that $I = \cap_{\cal Z} I_{\cal Z}\,$
(a more refined version of this is proved in
Theorem 6.2). If $I$ is generated by binomials,
then by Corollary 1.7 the ideal $ I_{\cal Z}$
is also generated by binomials.

The binomial ideals in  $k[{\cal Z}^\pm ]$
are completely classified by Theorem 2.1, and
Corollary 2.2 tells just when they are radical.  Thus
to classify all binomial algebraic sets $X$, it suffices
tell how the intersections of $X$ with the coordinate cells
can fit together.

\proclaim Theorem 4.1.
Let $k$ be any field.
An algebraic set $X \subseteq k^n\ $
is cut out by binomials if and only if the
following three conditions hold.
\item{(i)} For each coordinate cell
$(k^*)^{\cal Z}$, the algebraic
set $X\cap (k^*)^{\cal Z}$ is cut out by binomials.
\item{(ii)} The family of sets
$\, U \,= \,\{\,{\cal Z}\subseteq\,\{1,\ldots,n\} \ \ |\
\  X \, \cap \,(k^*)^{\cal Z} \, \not=\,
\emptyset \,\}\,$
is closed under taking intersections.
\item{(iii)} If
${\cal Z}, {\cal Z}' \in U$ and
${\cal Z}\subset {\cal Z}'$
then the coordinate projection
$(k^*)^{{\cal Z}'}\to (k^*)^{\cal Z}$
maps $\,X\cap (k^*)^{{\cal Z}'}\,$ onto a subset of
$\,X\cap (k^*)^{{\cal Z}}$.
\item{ }

We shall use the following definition and result.
A partially ordered set $U$ is a
{\it meet semilattice} if every finite subset
$\{ u_1,\ldots, u_m \} \subset  U$ has a unique greatest
lower bound in $U$. This lower bound is denoted
$u_1\wedge \dots \wedge u_m$ and
called the {\it meet} of $u_1,\ldots ,u_m$ in $U$.

\proclaim Lemma 4.2.
Let $U$ be a finite
meet semilattice and $R$ any commutative ring.
For each $u\in U$ let $J_u$
and $M_u$ be ideals in  $R$ such that
{a)} If $u\leq v$ then
$\sqrt{J_u}\subseteq \sqrt{J_v}$;  \
and
b)
$ \sqrt{M_{u \wedge v}} \subseteq \sqrt{M_u + M_v} $.
Under these assumptions, the two ideals
$$ \eqalign{
I_1 \,\, & = \,\, \bigcap_{u\in U}(J_u+M_u) \cr
\hbox{and} \qquad
I_2 \,\, & = \,\, \bigl( \bigcap_{u\in U}M_u \bigr) \,
+ \, \sum_{u\in U}\,
(J_u\cap\bigcap_{t\not\geq u}M_t) \cr}
$$
have the same radical
$\sqrt{I_1}=\sqrt{I_2}$.

\noindent {\sl Proof: }
To prove that $\sqrt{I_2} \subseteq \sqrt{I_1}$
it suffices to show that for all $u,v\in U$ we have
$$ J_u\cap\bigcap_{t\not\geq u}M_t
\quad \subseteq \quad J_v+M_v. $$
If $u\leq v$ then  $J_u \subseteq J_v$ by condition (a),
so $J_v$ contains the left hand side
and we are done.  If on the contrary
$u\not\leq v$, then $v$ is among the indices $t$ appearing
on the left hand side, so $M_v$ contains the left
hand side, and this suffices as well.

To prove that
$\sqrt{I_1} \subseteq \sqrt{I_2}$, choose a
 prime $P$ containing $I_2$. We must show that $P$
also contains $I_1$.  Let
$\,V \,= \,\{\,v \in U\ \ |\ \ M_v \subseteq P \}$.
{}From hypothesis (b) we see that
if $v,v'\in V$ then $v\wedge v' \in V$. Since
$P \supset \cap_{u\in U}M_u$, the set $V$ is nonempty.
Thus there
is a unique minimal element $w \in V$.
Since  $ P\supseteq I_2 \supseteq
 J_w\cap\bigcap_{t\not\geq w}M_t $
and $P$ does not contain any $M_t$ with $t\not\geq w$,
we see that $P$ contains $J_w$.  Thus
$P$ contains $J_w+M_w$, and with it $I_1$.\hfill\Box

\vskip .2cm

Here is the key part of the argument proving that
binomial ideals satisfy property (iii) of Theorem 4.1,
isolated for future use:

\proclaim Lemma 4.3.
Let  $\,
R :=  k[z_1,z_1^{-1},\dots,z_t,z_t^{-1}] \subset R' :=
k[z_1,z_1^{-1},\dots,z_t,z_t^{-1},y_1,\dots,y_s]\,$  be a
Laurent polynomial ring and a polynomial ring over it.
If  $B\subset R'$ is a binomial ideal
and $M \subset R'$ is a monomial ideal
such that $B+M$ is a proper ideal in $R'$, then
$$
(B+M)\,\cap\, R \,\,= \,\,B\,\cap \,R.
$$

\noindent {\sl Proof: }
Suppose $f\in (B+M)\cap R$.  The terms of
$f$ are invertible in $R'$. Since $B\!+ \!M \not= R'$,
no term of $f$ is in
$B \!+ \!M$.
Proposition 1.10 implies that $f \in B$.
\Box

\vskip .2cm

\noindent {\sl Proof of Theorem 4.1.}
Let $X \subset k^n$ be any algebraic set with
ideal $I \subset S$. Let  $U$ be the set of subsets ${\cal
Z} \subset \{1,\ldots,n\}$ such that $X\cap (k^*)^{\cal Z}$
is non-empty, or equivalently,  $I_{\cal Z} \not= S $.

Suppose $X$ is a binomial
algebraic set.  The ideal $I_{\cal Z}$
is binomial by Corollary 1.7, so
$X\cap (k^*)^{\cal Z}$ is cut out by binomials, proving
condition (i). To prove condition (ii) we must show
that if ${\cal Z}_1, {\cal Z}_2 \in U$ then
${\cal Z}_1\cap {\cal Z}_2 \in U$.
If on the contrary ${\cal Z}_1\cap {\cal Z}_2
\notin U$ then, for some integer $d$,
$$\,(\prod_{i\in {\cal Z}}x_i)^d\,\,\,\,
\in \,\,\,\, I \,+\,
M({\cal Z}_1 \cap {\cal Z}_2) \,\, = \,\,
I \, + \,M({\cal Z}_1) + M({\cal Z}_2). $$
Corollary 1.6 (b) implies that
$(\prod_{i\in {\cal Z}}x_i)^d $
is in either in $ \,I + M({\cal Z}_1)\, $
or in $\,I + M({\cal Z}_2) $.
Consequently either  $I_{{\cal Z}_1}$
or $I_{{\cal Z}_2}$ is the unit ideal in $S$,
contradicting our assumption.

Write $k[{\cal Z}]$ for the polynomial ring
$k[\{x_i\}_{i \in {\cal Z}}]$.
The algebraic form of condition (iii) is the statement
that if ${\cal Z},{\cal Z}'\in U$
with ${\cal Z} \subset {\cal Z}'$ then
$I_{\cal Z}\cap k[{\cal Z}]
\subseteq I_{{\cal Z}'}$.
Since
$I_{{\cal Z}'} =
(I_{{\cal Z}'}:\prod_{i\in {\cal Z}} x_i)$,
it suffices to
prove this condition after inverting
the $x_i$ for $i\in {\cal Z}$.  That is, if we set
$R' = k[{\cal Z}^\pm][\{x_i\}_{i\notin {\cal Z}}]$,
then we must show that
$$ (I+M({\cal Z}))R' \,\cap \,k[{\cal Z}^\pm]
\quad \subseteq \quad I_{{\cal Z}'}R'. $$
Since ${\cal Z} \in U$, the ideal
$
(I+M({\cal Z}))R'
$
is proper, and we may apply Lemma 4.3 to get
$(I+M({\cal Z}))R'\cap k[{\cal Z}^\pm]=
IR'\cap k[{\cal Z}^\pm].
$
Since $I \subseteq I_{{\cal Z}'}$, we are done.

Conversely, suppose that $X$ is any algebraic set
satisfying conditions (i),(ii) and (iii).
We must show that the ideal $I$ of $X$ is generated
by binomials. We have already remarked that
$I = \cap_{{\cal Z}\in U} I_{\cal Z}$.
Note that $U$ is a partially ordered set under the
inclusion relation for subsets of  $\{1,\dots,n\}$. By
condition (ii) the set $U$ is closed under intersections,
so $U$ is a meet semilattice.
For ${\cal Z}\in U$ we set  $ \,J( {\cal Z})\, :=
\,( I_{\cal Z}\cap k[{\cal Z}])  S \,$
and, as before,  $M ({\cal Z})
\, = \, (\{x_i\ |\ \ i \notin {\cal Z}\})$.
We shall apply Lemma 4.2 to these ideals.
Hypothesis (b)  of Lemma 4.2 is obvious from
the definition of $M ({\cal Z})$, and hypothesis
(a) is implied by the algebraic form of condition (iii)
given above. The ideal
$I_1 = \bigcap_{{\cal Z}\in U}(J({\cal Z})+M({\cal Z}))$
equals $\bigcap_{{\cal Z}\in U} I_{\cal Z} = I$.
Each $J({\cal Z})$ is a binomial ideal by Corollary 1.3,
and each $M({\cal Z})$ is a monomial ideal.
Hence each term in the sum
$$ I_2 \quad = \quad \sum_{{\cal Z}\in U}
\, \bigl(\,J({\cal Z})\,\,\cap \,\bigcap_{{\cal Z}'
\not\supseteq {\cal Z}} M ({\cal Z}') \,\bigr) $$
is a binomial ideal by Corollary 1.5.
This shows that $I_2$ is binomial. Theorem 3.1
now implies that $\sqrt{I_2} = \sqrt{I_1} = I$
is binomial, as claimed.
\hfill \Box

\proclaim Problem 4.4. \ (Find the generators) \ \rm
\ \ \ In the application of Lemma 4.2 made
in the proof of Theorem 4.1, are the
ideals $I_1$ and $I_2$ actually equal?
This is the case when the set $U$ is totally ordered
and in other examples we have tried,
such as the following:

\proclaim Example 4.5. \ (Subsets of the
vertex set of the coordinate cube) \hfill \break \rm
For each ${\cal Z}\subseteq\{1,\ldots,n\}$ let $p_{\cal Z}$ be the
point whose $i^{\rm th}$ coordinate is $1$ if $i\in {\cal Z}$ and
$0$ otherwise. Let $U$ be a collection of subsets of
$\{1,\ldots,n \}$. The finite algebraic set
$$ X_U \quad := \quad \{\,p_{\cal Z}\ \ |\ \ {\cal Z} \in
 U \, \} \quad \subset \quad k^n $$
is cut out by binomials if and only if $U$ is closed under
taking intersections.
We remark that for {\sl any} collection
$U$ of subsets, the ideal of $X_U$ is
generated by the $n$ binomials
$x_i(x_i-1)$ for $1\leq i\leq n$ (these generate the ideal
of all the $2^n$ points $p_Z$)
and the $card(U)$ elements
$$
\prod_{i\in {\cal Z}}(x_i-1)
\prod_{i\notin {\cal Z}}x_i \qquad
\hbox{for \ $Z \in U$.} $$
\vskip .2cm

\proclaim Example 4.6.
A binomial algebraic set whose top-dimensional part
is not binomial. \break \rm
Consider the following three binomial
varieties in affine $4$-space $k^4$:
$$ \eqalign{  & V_1 \,= \,V(x_1x_2-1,x_3,x_4), \quad
\hbox{a hyperbola in the cell} \,\, (k^*)^{\{1,2\}}; \cr
& V_2 \,= \,V(x_1,x_2,x_3x_4-1), \quad
\hbox{a hyperbola in the cell}\,\, (k^*)^{\{3,4\}}; \cr
& V_3 \,= \,V(x_1,x_2,x_3,x_4) ,\quad
\hbox{ the unique point in the cell} \,\,
(k^*)^\emptyset .\cr}$$
The union of these varieties is defined by the
binomial ideal
$$ \eqalign{ &
I(V_1\cup V_2\cup V_3) \quad = \quad
I(V_1)\,\cap \,I(V_2)\,\cap \,I(V_3) \quad = \cr
&= \quad  \bigl( \,x_1^2x_2-x_1,\, x_1x_2^2-x_2,
\,x_3^2x_4-x_3, \,x_3x_4^2-x_4, \,x_1x_3,\ x_1x_4,\
x_2x_3,\ x_2x_4 \,\bigr) .\cr}$$
However, the union of $V_1$ and $V_2$, the top-dimensional
components, is not cut out by  binomials. Its ideal
$I(V_1\cup V_2)$ has the reduced Gr\"obner basis
$$\bigl\{ x_1x_2+x_3x_4-1,\, x_3^2x_4-x_3 ,
\,x_3x_4^2-x_4 , \,x_1x_3,\ x_1x_4,\ x_2x_3,\ x_2x_4
\, \bigr\} .$$
By homogenizing these equations we get a projective
binomial scheme with the same property. Note also that
$\,(I(V_1\cup V_2\cup V_3) :(x_1,x_4)) = I(V_1\cup V_2) $,
so this ideal also exhibits the phenomenon of
Example 1.8. \Box

\vskip .2cm

\proclaim Example 4.7. \
Face rings of polyhedral complexes.
\rm (cf.~Stanley [1987], \S 4) \hfill \break
By a {\it lattice polytope}
in ${\bf R}^m$ we mean the convex hull of a finite
subset of ${\bf Z}^m$. A (finite, integral)
{\it polyhedral complex} $\Delta$  is a finite set of
lattice polytopes in ${\bf R}^m$, satisfying
\item{(i)} any face of a polytope in
$\Delta$  is a polytope in $\Delta$;
\item{(ii)} any two of the polytopes
in $\Delta$ intersect in a set that is a face
of each of them.

The polytopes in $\Delta$ are called {\it faces} of
$\Delta $.  The maximal faces are called {\it facets}.
We write ${\cal F}(\Delta)$ for the set of facets of
$\Delta$. For each face $P \in \Delta$ we define a cone
$$
C_P \, =\, \{(a_1,\dots,a_m,b) \in {\bf R}^{m+1} \ \  | \
\  (a_1,\dots,a_m,b)=(0,\dots,0)  \ {\rm or}\
(a_1/b,\dots,a_m/b)\in P \}.
$$
Stanley [1987] defines the {\it face ring}
$k[\Delta]$ of $\Delta$ to be the ring having vector
space basis over $k$ the set of monomials
$ \{ y^\alpha \ |\  \alpha \in C_P \cap {\bf Z}^{m+1}
\,\,{\rm for\  some\ } \,
P \in \Delta \} $ with multiplication
$$ y^\alpha y^\beta \quad = \quad
\cases{y^{\alpha + \beta }, &if
$\alpha , \beta  \in C_P$ {\rm
for\  some\ }  $\, P \in \Delta $;\cr
0, &otherwise.}
$$
If $\Delta $ has a single facet $P$,
then the face ring of $\Delta$ is the homogeneous
coordinate ring $k[P]$ of the projective toric
variety associated with the lattice polytope $P$
(see for example Fulton [1993] or Sturmfels [1991]).
We may represent it as
$$
k[P] \quad = \quad  k[\{x_{i}\}_{i\in G(P)}]/I(P)
$$
where the $x_{i}$ are variables
indexed by the unique minimal set $G(P)
\subset {\bf Z}^{m+1}$ of additive  generators for the
semigroup  $\,C_P \,\cap \,{\bf Z}^{m+1} \,$
and $I(P)$ is the binomial prime ideal of relations
among the monomials $y^\beta$ for $\beta\in G(P)$.

More generally, let $\,G(\Delta ) := \cup_{P \in {\cal F}(\Delta)} G(P) $.
We may represent the face ring of $\Delta$ as
$$ k[\Delta] \quad = \quad k[\{ x_i \} _{i\in
G(\Delta)}]/I(\Delta)$$ for some ideal $I(\Delta)$. This
ideal is an intersection of binomial primes satisfying
Theorem 4.1, so it is generated by binomials. The
following more precise result is implicit in Stanley
[1987];  the proof was communicated to us privately by
Stanley. Its geometric interpretation is
that the projective scheme $ Proj(k[\Delta]) $ is the
reduced union of the toric varieties $Proj(k[P])$,
glued along orbit closures corresponding
to intersections of facets in $\Delta$.

\proclaim {Proposition 4.8}.
The ideal $I(\Delta)$ defining the face ring
$k[\Delta]$ is the intersection of the binomial primes
$\,  I(P)+ (\{x_{i}\}_{i\in G(\Delta)\setminus G(P)})$,
where $P$ ranges over the set of facets ${\cal F}(\Delta)$.
The ideal $I(\Delta)$ is generated by  $\sum_{P\in {\cal F}(\Delta)}I(P)$
together with all the monomials  $x_{i_1}\cdots x_{i_s}$ such that
$i_1,\dots ,i_s$ do not all lie in any facet of $\Delta$.

\noindent {\sl Proof}. The $k$-basis given for $k[\Delta]$
in the definition is a subset of the natural
vector space basis of
$\prod_{P\in{\cal F}(\Delta)} k[P]. $
The description of the multiplication gives
an inclusion of  $k$-algebras
$\, k[\Delta] \subset
\prod_{P\in{\cal F}(\Delta)} k[P] $.
The ideal $I(\Delta)$ is by definition
the kernel of the natural map
$k[\{ x_i \} _{i\in G(\Delta)}]
\to \prod_{P\in{\cal F}(\Delta)} k[P] $.
It follows that  $I(\Delta)$ is the intersection
of the ideals
$ J(P) := ker \bigl(\ k[\{ x_i \} _{i\in G(\Delta)}]
\to k[P]\ \bigr) $ for $P\in {\cal F}(\Delta)$,
and it is immediate that $J(P) =
I(P)+ (\{x_{i}\}_{i\in G(\Delta)\setminus G(P)})$.
This proves the first assertion.

Let $I$ be the ideal generated by
$\sum_{P \in {\cal F}(\Delta) }I(P)$
and the non-facial monomials  $x_{i_1}\cdots x_{i_s}$.
The inclusion $I \subseteq I(\Delta)$ is evident,
so we get a  surjection from
$R := k[\{ x_i \}_{i\in G(\Delta)}]/I$
onto $k[\Delta]$. Each non-zero monomial in $R$ is mapped
to a monomial $y^\beta $ in $k[P]$ for some $P$. Any two preimages
of $y^\beta$ differ by an element of $I(P) \subset I $,
hence they are equal in $R$. This shows that
the surjection is injective as well,
and therefore $I$ equals $I(\Delta)$, as desired.
\Box

\vskip .1cm

\proclaim Problem 4.9.  \ Intersections of binomial
ideals. \ \rm \hfill \break
It would be nice to have a result like Theorem 4.1 for
the  intersections of arbitrary binomial ideals,
not just radical binomial ideals.
A first step might be to answer the following
question: Which sets of primes
can be the set of associated primes
of a binomial ideal~?

In some cases a fairly straightforward
generalization to schemes of Theorem 4.1
seems to be all that is necessary.  For example,
the following union of three lines, contained
in the closures of the
$\{x_3\}$-cell, the $\{x_2,x_3\}$-cell, and the
$\{x_1,x_3\}$-cell respectively, is binomial:
$$
(x_1,x_2)\,\cap \,(x_1,x_2-x_3)\,\cap \,(x_2,x_1-x_3).
$$
If we thicken the line in the $\{x_2,x_3\}$-cell then
we get a scheme that is not binomial:
$$
(x_1,x_2)\,\cap\, (x_1^2,x_2-x_3)\,\cap \,(x_2,x_1-x_3)
$$
However, if we also thicken the line in the $\{x_3\}$-cell
enough so that the line in the
${x_2,x_3}$-cell projects into it,
$$
(x_1^2,x_2)\,\cap \,(x_1^2,x_2-x_3)\,\cap\, (x_2,x_1-x_3)
$$
then again we get a binomial scheme.\hfill \Box

\vskip .4cm

\beginsection 5. Some binomial ideal quotients

The theory of binomial ideals would be much easier if the
quotient of a binomial ideal by a binomial were again a binomial
ideal.  Here is a simple example where this fails:

\proclaim {Example 5.1}. \rm Let $ \, I \,
= (x_1 - y x_2, x_2 - y x_3, x_3 - y x_1) \,
\subset \, k[x_1,x_2,x_3,y]$. The ideal
$\, (I: (1-y)) \,=\,
(x_1 + x_2 + x_3,\,
x_2^2 + x_2 x_3 + x_3^2, \,
x_2 y + x_2 + x_3, x_3 y - x_2)\,$ is not binomial:
the given generators form a reduced Gr\"obner basis,
so Corollary 1.2 applies.

By reducing problems to coordinate cells
$(k^*)^{\cal Z}$ as in Section 4, we
can often assume that some variables
are nonzerodivisors modulo a given
binomial ideal $I$.  In such a case certain ideal quotients
of $I$ by a binomial are again binomial.  The results we
are about to  prove are the main technical tools for our
study of primary decompositions.

The ordinary powers of a binomial are not binomials.
However, there is a natural binomial operation that has
many features in common with taking powers:
If $b = x^{\alpha} -  c x^{\beta}$ is a binomial, and
$d$ is a positive integer, then
we set   $\, b^{[d]}\,:= \,x^{d \alpha} -  c^d x^{d
\beta}\,$  and call it the $d^{\rm th}$
{\it quasi-power} of $b$. Note that if $d|e$ then
$b^{[d]}\,| \, b^{[e]}$.

\proclaim {Theorem 5.2}. Let $I $  be a binomial
ideal in $S = k[x_1,\ldots,x_n]$
and $<$ a monomial order on $S$.
Suppose $b:=x^\alpha-ax^\beta$ is a binomial  and $f \in S$
such that $bf \in I$ but
$x^\alpha$ is a nonzerodivisor mod $I$.
Let $f_1+\dots+f_s$ be the normal form of $f$
modulo $I$ with respect to $<$.
If $d$ is a sufficiently divisible
positive integer, then \item{(a)} the binomials
$\,b^{[d]}f_j \,$ lie in $\, I \,$ for
$\, j=1,\ldots,s $.
\item{(b)} $ (I:b^{[d]}) $
is generated by monomials modulo $I$, and is thus a binomial ideal.
\item{(c)}  Let $p = char(k)$. If $p = 0$ let $q=1$, while
if $p > 0$ let $q$ be the largest power of $p$
that divides $d$. If $e$ is a divisor of $d$
that is divisible by $q$,
then $ (I:(b^{[d]}/b^{[e]})) $ is a binomial ideal.

\noindent {\sl Proof.  }
(a): \ To say that $f_1+\dots+f_s$ is the normal form of $f$
modulo $I$ means that $f \equiv f_1+\dots+f_s\ (mod\ I)$
and that the $f_i$ are terms
not in $in_<(I)$.
By Proposition 1.1 (b), the normal form of each term
$x^\alpha f_i$ or $ax^\beta f_j$  is a term.
Since $x^\alpha$ is a nonzerodivisor modulo $I$,
the terms $x^\alpha f_i$, $i=1,\ldots,s$,
have distinct non-zero normal forms.
The equation $ \, (x^\alpha-ax^\beta)f \equiv 0 \,\,mod\,\, I\,$
shows that for each of the $s$ terms $x^\alpha f_i$
there exists a term $ax^\beta f_j$ with the same normal form,
and by counting we see that $j$ is uniquely determined.
Since $x^\alpha$ is a nonzerodivisor modulo $I$ we may
invert it and set $r := ax^\beta /x^\alpha$.  The discussion above
amounts to saying that (modulo $I$) multiplication by $r$ induces
a permutation of the terms $f_i$ of $f$.  We let $d$ be any
positive integer divisible by the order of this permutation.
Multiplication by $r^d$ induces the identity permutation,
so $(1-r^d)f_i \in I$.  Multiplying by $x^{d\alpha}$ we get part (a).

\vskip .1cm

(b): \ If $d$ and $d'$ are positive integers such that
$d$ divides $d'$, then the binomial
$ b^{[d]}$ divides the binomial
$ b^{[d']}$, so that $\,\bigl( \,I \,:\,b^{[d]} \,\bigr)
\, \subseteq \,\bigl( \,I \,:\,b^{[d']} \,\bigr)$.
Since $S$ is Noetherian we may
choose $d$ sufficiently divisible so that
equality holds for all integers $d'$
divisible by $ d $. We claim that for such a choice
of $d$ the conclusion of part (b) is  satisfied.
Let $f\in (I:b^{[d]})$.
By induction on the number of terms $f_i$ in the normal
form of $f$ modulo $I$, it suffices to show that the
first term $f_1$ is
in $(I:b^{[d]})$.
By part (a) applied to $b^{[d]}$,
there is an integer $d'$ such that
$f_1\in (I:b^{[dd']})$.
By the choice of $d$ we have $f_1\in (I:b^{[d]})$
as desired.

\vskip .1cm

For the proof of part (c) we use a general fact:

\proclaim Lemma 5.3. Let R be any commutative ring and
$f,g\in R$. If $(f,g) = R$ then $(0:g) = (0:fg)f$.

\noindent {\sl Proof}. It is immediate that $(0:g)\supseteq (0:fg)f$.
For the opposite inclusion, suppose $x \in (0:g)$.  Since $(f,g) = R$
we may write $1=af+bg$ with $a,b\in R$, so we have
$x = xaf+bgx = xaf$.  Since
$xa\cdot fg = xg\cdot af = 0\cdot af = 0$ we get $x = xaf \in (0:fg)f$.
\Box

\vskip .2cm

\noindent {\sl Proof of Theorem 5.2 (c):  }
We apply Lemma 5.3 to the ring
$\,R = (S/I)[1/x^\alpha]\,$ with
$f = b^{[e]}$
and $g = b^{[d]}/b^{[e]}$. The hypothesis $(f,g)=R$
of Lemma 5.3 holds because over the algebraic
closure $\bar k$ of $k$ we have the factorizations
$$ \eqalign{
f \quad &=\quad
\prod_{\eta\in k^*,\, \eta^{e/q}=1}
\!\!\!\!\!\! (x^{q\alpha}-\eta a^qx^{q\beta})\cr
g \quad &=\quad
\prod_{\zeta\in k^*,\, \zeta^{d/q}=1,\,\, \zeta^{e/q}\neq 1}
\!\!\!\!\!\!\!\! (x^{q\alpha}-\zeta a^qx^{q\beta}).
}$$
Any factor of $f$ together with any factor
of $g$ generates the unit ideal, and hence
$ (f, g) = R. $

If $J$ is an arbitrary ideal of $S$ then the
preimage of $JR$ in $S$ is $((I+J):(x^\alpha)^\infty)$.
If $J = (I:g)$ then $I\subseteq (I:g)$. Since
$x^\alpha$ is a nonzerodivisor modulo $I$ it is
also a nonzerodivisor modulo $(I:g)$. Thus the preceding
formula simplifies, and the preimage of
$(I:g)R$ in $S$ is equal to $(I:g)$. Applying Lemma 5.3 and
pulling everything back to $S$, we get $$
(I:g) \,\,=\,\,  \bigl( \,(I+(I:fg)f)\,:\,(x^\alpha)^\infty \bigr)
\qquad \hbox{in $\,\,S$.} $$
By part (b), the ideal $(I:fg) = (I : b^{[d]})$
is generated modulo $I$ by  monomials.
Since $f = b^{[e]}$ is a binomial, $I+(I:fg)f$ is binomial.
By Corollary 1.7,
the quotient $\, ((I+(I:fg)f):(x^\alpha)^\infty)\, $ is a
binomial ideal, and thus $(I:g)$ is binomial as  desired.
\Box

\proclaim {Example 5.1, continued}. \rm For $ \, I \,
= (x_1 - y x_2, x_2 - y x_3, x_3 - y x_1) $, the
ideals $ (I:(1-y^3))\,=\, ( x_1,x_2,x_3 )
\,$ and $\,(I : (1+y+y^2)) = (x_1-x_3,\ \ x_2-x_3,\ \
x_3 y-x_3) \,$ are binomial.

\vskip .2cm

\noindent {\bf Example 5.4.} \ \
The hypothesis that $x^\alpha$ is a non-zerodivisor
is necessary for Theorem 5.2~(b) to hold. For instance,
consider the radical binomial ideal
$$ \eqalign{ I \quad = \quad & (u x-u y, u z-vx, vy-vz) \cr
= \quad  & (x, y, z)\,\cap\,
(u, v)\,\cap\,
(u, x, y-z)\, \cap \,
(v, z, x-y) \, \cap \,
(x-y, y-z, u-v )
}
$$
in $k[x,y,z,u,v]$.
Both $u$ and $v$ are zerodivisors mod $I$.
For each positive integer $d$ we have
$$ \bigl( I: (u^d-v^d)\bigr)\quad
= \quad (x-y+z,\, u z,\, yz-z^2, \,vy-vz ).$$
This quotient is not a binomial ideal.
 \Box

\vskip .1cm

The following two results on quasi-powers will be
used in the proof of Theorem 7.1.

\proclaim Proposition 5.5.
Let $I \subset S$ be a binomial ideal,
and let $b = x^{\alpha} -  a x^{\beta}$ be a binomial such
that $ x^{\alpha}$ is a nonzerodivisor modulo $I$. For
sufficiently divisible positive integers $d$ we have $$
\bigl(\, I \,:\, b^{[d]} \,\bigr) \quad = \quad \bigr(\, I
\, : \, (b^{[d]})^2 \,\bigr) .$$

\noindent {\sl Proof: }
By Theorem 5.2 (b), the quotient
$(I:b^{[d]})$ is generated by monomials mod $I$,
and this ideal is independent of $d$ for sufficiently
divisible $d$.  Using Theorem 5.2 (b) again we see that
$(I:(b^{[d]})^2) = ((I:b^{[d]})\, :\, b^{[d]})$
is generated by monomials mod $I$, and it suffices
to show that if $m \in (I:(b^{[d]})^2)$ is a monomial
then $m \in (I:b^{[d]})$.  By Proposition
1.1 (b), the normal form of $m$
mod $I$ is a term, and we may assume that
it equals $m$. Now $b^{[d]}$ annihilates $b^{[d]}m$ mod $I$
by hypothesis, so Theorem 5.2 (a) implies that $\,b^{[d]}$
annihilates $x^{d \alpha}m$.  Since
$x^{d \alpha}$ is a nonzerodivisor mod $I$, we see that
$b^{[d]}$ annihilates $m$ mod $I$.
\Box

\proclaim Corollary 5.6.
Let $I$ be a binomial ideal in $S$, and let
${\cal Z} \subseteq \{1,\dots,n\}$ be a subset such that
$x_i$ is a nonzerodivisor modulo $I$ for each $i \in {\cal Z}$.
If $\sigma$ is a partial character on ${\bf Z}^{\cal Z}$
and
$\sigma_d$ is the restriction of $\sigma$ to $\,d L_\sigma$,
then for sufficiently divisible integer $d$ we have
$$ \bigl( I : I_+ (\sigma_d) \bigr) \quad = \quad
   \bigl( I : I_+ (\sigma_d)^\infty \bigr) .$$

\noindent {\sl Proof:  }
Consider the ideal  $I_+ (\sigma_d)$
in $k[{\cal Z}]$.
It is generated by the  $d^{\rm th}$ quasipowers
of all binomials in $I_+ (\sigma)$, and of course
a finite set $\{b_1,\ldots,b_s \} \subset
I_+(\sigma)$ suffices.
 For each $i$ the two monomials of $b_i$ are
nonzerodivisors mod $I$ because they are monomials
in $k[{\cal Z}]$.
By Proposition 5.5 we have  $\, ( I \,:\, b_i^{[d]} ) \, = \,
     ( I \, : \, (b_i^{[d]})^2 ) $.
Since the quasipowers $b_i^{[d]}$ generate
$I_+(\sigma_d)$, we get
$$ \bigl( I : I_+ (\sigma_d) \bigr) \quad = \quad
   \bigcap_{i=1}^s \bigl( I \,:\, b_i^{[d]} \bigr) \quad = \quad
   \bigcap_{i=1}^s \bigl( I \, : \, (b_i^{[d]})^2 \bigr) \quad
    \supseteq \quad  \bigl( I : I_+ (\sigma_d)^2 \bigr) .$$
The reverse inclusion is obvious, and
it implies the desired result.
\Box

\vskip .6cm

\beginsection 6. Associated primes, isolated components
and cellular decomposition

The decompositions in a univariate polynomial ring
$$\eqalign{
(x^d-1) \quad &= \quad (x-1)\cap (x^{d-1}+\dots+x+1)\cr
(x^{d-1}+\dots+x+1) \quad &= \quad
\bigcap_{\zeta^d=1,\ \zeta\neq 1} (x-\zeta)
\cr} \eqno (6.1) $$
show that in order for the associated primes of
a binomial ideal to be binomial we
must work over a field $k$ containing the roots of
unity.  Further, for the
minimal primes of $(x^d-a)$ to have the  form
given in Corollary 2.4, the scalar $a\in k$ must
have all its $d^{\rm th}$ roots in $k$.
This is the reason why
$k$ is taken to be algebraically closed
in the following theorem.

\proclaim Theorem 6.1. Let $k$ be an algebraically
closed field. If $I$ is a binomial ideal in
$ S = k[x_1,\dots,x_n]$, then every associated prime
of $I$ is generated by binomials.

\noindent {\sl Proof:  }  If
$\, I \, = \, I_+(\rho) \, = \,
I(\rho) \cap k[x_1,\dots,x_n] \,$
for some partial character $\rho$ on ${\bf Z}^n$
then Corollary 2.3 implies the desired result.
We may therefore assume that there is
a variable $x_i$ such that $(I:x_i) \neq I$.
If $x_i \in I$ we may reduce modulo $x_i$
and do induction on the
number of variables. Hence may assume that
$x_i \notin I$. From the short exact sequence
$$ 0 \,\, \rightarrow \,\,
S/(I:x_i) \,\, \rightarrow \,\,
S/I \,\,\rightarrow\,\, S/(I,x_i)\,\, \rightarrow \,\,
0 \eqno (6.2) $$
we see that
$Ass(S/I) \subseteq Ass(S/(I:x_i))\cup Ass(S/(I,x_i))$.
By Noetherian induction and Corollary 1.7,
both of these sets
consist of binomial primes.
\Box

\vskip .1cm

Corollary 2.3 does primary decomposition
for binomial ideals whose
associated points are all contained in the open cell away
from the coordinate hyperplanes.  This suggests dividing
up the primary components  according to
which coordinate cells they lie in.
We define an ideal $I$ of $S$ to be {\it cellular}
if, for some  ${\cal Z}\subseteq \{1,\dots,n\}$, we have
$I = (I:(\prod_{i\in {\cal Z}} x_i)^\infty)$
and $I$ contains a power of
$\, M({\cal Z}) = ( \{x_i\}_{ i \not\in {\cal Z}} )$.
This means that the scheme defined
by $I$ has each of its associated points in the cell
$(k^*)^{\cal Z}$.

Given any ideal
$I\subseteq S $ we can manufacture
cellular ideals from $I$ as follows. For each
vector of positive integers $d = (d_1,\ldots,d_n) $
and each subset ${\cal Z} $ of $\{1,2,\ldots,n\}$ we set
$$ I_{\cal Z}^{(d)} \quad := \quad
\biggl( \,(\,I\,+\,(\{x_i^{d_i}\}_{i \notin {\cal Z}})) \, :\,
(\prod_{j \in {\cal Z}} x_j)^\infty \,\biggr). \eqno (6.3)
$$
For $d=(1,\ldots,1)$ we have
$I_{\cal Z}^{(d)} = I_{\cal Z}$, the ideal
considered for $I$ radical in Section 4.

\proclaim {Theorem 6.2}. The ideal $I_{\cal Z}^{(d)}$
is a cellular binomial ideal for all $I$, $d$ and
${\cal Z}$.  For distinct ${\cal Z}$ and ${\cal Z}'$
the sets  of associated primes
$Ass(I_{\cal Z}^{(d)})$ and $Ass(I_{{\cal Z}'}^{(d)})$
are disjoint. If the integers $d_i$ are chosen sufficiently
large, then
$$ I \quad =  \, \bigcap_{{\cal Z}\subseteq \{1,\ldots,n\}}
\! \! I_{\cal Z}^{(d)}. \eqno (6.4) $$
Thus an irredundant primary decomposition
of $I$ is obtained from given primary decompositions
of the $I_{\cal Z}^{(d)}$ by deleting redundant components.
Equation (6.4) holds, in particular, if for some  primary
decomposition $I = \cap Q_j$ we have  $ x_i \in \sqrt{Q_j
}$ if and only if $ x_i^{d_i} \in Q_j $ for all $i$ and
$j$.

We say that the binomial ideals in (6.3) form a
{\it cellular decomposition} of  $I$.

\vskip .1cm

\noindent {\sl Proof}.
The $I_{\cal Z}^{(d)}$ are binomial by Corollary 1.7.
They are obviously cellular.  The
primes associated to $I_{\cal Z}^{(d)}$  contain the variable
$x_i$ if and only if $i\in {\cal Z}$, and this shows that
the sets of associated primes $Ass(I_{\cal Z}^{(d)})$ are
pairwise disjoint.

We next show that if the $d_i$ are chosen
to have the property specified with respect
to a primary decomposition $I = \cap Q_j$,
then $I$ is the intersection of the
ideals $I_{\cal Z}^{(d)}$. Our assertion about
primary decomposition follows at once from this.
Since $I$ is obviously contained in the
intersection of the $I_{\cal Z}^{(d)}$,
it suffices to prove that for each
$f \in S \setminus \! I$, there exists an index
set ${\cal Z} \subseteq \{1,\ldots,n\}$
such that $f\notin I_{\cal Z}^{(d)}$.

Let $m = x_{i_1} x_{i_2}  \cdots x_{i_r}$ be a maximal
product of variables such that
$\,f \not\in (I:m^\infty)\,$ and define
${\cal Z} := \{i_1,\ldots,i_r\}$.
We have $\,(I:m^\infty) \, = \,\cap \,( Q_j :m^\infty)$.
Thus there exists a primary component $Q_s$ with
$ f \not\in ( Q_s :m^\infty)$. It follows
that $ ( Q_s :m^\infty) \not= S$, hence
$\, ( Q_s :m^\infty) = Q_s \,$ and $\,f \not\in Q_s$.

By the maximality in our choice of $m$,
each variable $x_j$ with $j \notin {\cal Z}$ has a power
throwing $f$ into $\,(I:m^\infty) \,$ and hence
throwing $f$ into $\,(Q_s :m^\infty) \, = \,Q_s $.
We see that the variables $x_j$, $j \notin {\cal Z}$, are
zero-divisors modulo $Q_s$, hence they are
nilpotent modulo $Q_s$. This implies $x_j^{d_j} \in Q_s$
for $j \notin {\cal Z}$. This proves that
$$
    Q_s \quad = \quad
 ( Q_s : m^\infty ) \quad = \quad
\bigl( \,(Q_s\,+\,(\{x_j^{d_j}\}_{j \notin {\cal Z}})) \, :\,
(\prod_{j \in {\cal Z}} x_j)^\infty \,\bigr) .
$$
This ideal contains  $I_{\cal Z}^{(d)}$, as can be seen
from (6.3), and therefore $f \not\in I_{\cal Z}^{(d)}$.
\Box

\vskip .2cm

\noindent {\bf Problem 6.3. }
It would be nice to have a criterion for when
the $d_i$ are large enough for (6.4)
that does not require the knowledge of a primary decomposition
$\,I \, = \, \cap Q_j$. Perhaps such a criterion can
be found using the methods in the proof in the
effective Nullstellensatz given by Koll\'ar [1988].
We remark that the conditions
$(I : x_i^{d_i}) = (I:x_i^\infty)$
are not sufficient. For instance, let
$\, I \, := \, (
 x_1 x_4^2 - x_2 x_5^2 ,\,
x_1^3 x_3^3 - x_2^4 x_4^2 ,\,
x_2 x_4^8 - x_3^3 x_5^ 6 )\,$
and $\, d = (2,2,0,4,5)$. Then
$\, ( I :  x_i^{d_i}) = ( I :  x_i^\infty) $
for all $i$, but $I$ is properly contained in
$\,\cap_{\cal Z} I_{\cal Z}^{(d)} $.
(There are only two cellular components in this example:
$\,{\cal Z} = \{1,2,3,4,5\} $ and ${\cal Z} = \{3 \} $).

\vskip .3cm

The main results of this section are the following
theorem and its corollary, which say that in certain cases
the localization of a cellular binomial ideal is binomial.
If $I,J$ are ideals of $S$, then we define
$I_{(J)}$ to be the intersection of all those
primary components of $I$ that are contained in
some minimal prime of $J$.  (The notation is motivated
by the fact that if $J$ is prime then
$I_{(J)}= S\cap IS_J$, where $S_J$ is the usual
localization.)

\proclaim {Theorem 6.4}. If $I$ and $J$ are
binomial ideals in $ S=  k[x_1,\dots,x_n]$ that
are cellular with respect to the same index set
${\cal Z} \subseteq\{1,2,\ldots,n\}$,
 then the ideal $I_{(J)}$ is  binomial.

\noindent {\sl Proof}. We may harmlessly replace
$J$ by its radical and thus assume that
$\,J\,=\, M({\cal Z})\,+\,I_+(\sigma)$
for some partial character $\sigma$ on ${\bf Z}^{\cal Z}$.
By Corollary 1.2 we may assume that $k$ is
algebraically closed. Further, by Noetherian induction,
we may suppose that the result is
true for any binomial ideal strictly containing $I$.

If all the associated primes of $I$
are contained in a minimal prime of
$J$, then $I= I_{(J)}$ and we are
done.  Else let
$\,P \,=\, M({\cal Z})\,+\,I_+(\rho) \,$
be a prime associated to $I$ that is not
contained in any minimal prime of $J$.
We consider the following sublattice of
${\bf Z}^{\cal Z}$,
$$ L \quad := \quad \{\,m\in L_{\sigma}\cap L_{\rho} \,\,:
\,\,\sigma(m)=\rho(m)\,\}, \eqno (6.5) $$
and we distinguish two cases:

{\sl Case 1:  $L$ has finite index in $L_{\rho}$.}
Since $L\subseteq L_\sigma$
we see in this case that $L_\rho \subseteq Sat(L_\sigma)$.
We first claim that  $L \neq L_\rho \cap L_\sigma$.
In the contrary case we could define
a partial character $\tau$ on $L_\rho + L_\sigma$ by the
formula $\tau(m+ \tilde m) = \rho(m)+\sigma(\tilde m)$ for
$m\in L_\rho$ and $\tilde m\in L_\sigma$.  Since
$k^*$ is a divisible group, one of the
saturations $\sigma'$ of $\sigma$ would extend $\tau$, and
thus $I_+(\rho)$ would be contained in the minimal prime
$I_+(\sigma')$ of $I_+(\sigma)$, contradicting our
hypothesis and establishing
the claim.  It follows that we may choose an element
$m \in L_\rho \cap L_\sigma$ that is not in $L$,
so that $\sigma(m)\neq \rho(m)$.
The binomial $b:=x^{m_+}-\sigma(m)x^{m_-}$ is in $J$ but
not in $P$.

Since the index
of $L$ in $L_{\rho}$ is finite, there is a  root of unity
$\zeta$ such that $\rho(m)=\zeta \sigma(m)$.  If $d$ is a
sufficiently divisible integer,
 and $q$ is the largest power
of the characteristic of $k$ that divides $d$
\ (or $q= 1$ if $char(k)=0$),
then the ratio of quasi-powers
$ g \,\, = \,\, b^{[d]}/b^{[q]} $
lies in $P$ but not in any minimal prime of $J$.
By Theorem 5.2 (c), the ideal $I':=(I:g)$ is
binomial.  It is larger than $I$ because
$g \in P \in Ass(S/I)$. On the other hand,
$I'_{(J)}=I_{(J)}$ because $g$ is not
in any minimal prime of $J$,
so we are done by Noetherian induction.

\vskip .1cm

{\sl Case 2:  $L$ does not have finite index in $L_{\rho}$.}
We may choose an element $m \in L_{\rho}$ whose image in
$L_{\rho}/L$ has infinite order.
Set $b = x^{m_+} - \rho(m) x^{m_-}$. For any integer
$d > 0 $, the quasi-power $b^{[d]}$ is in
$P$ but not in any minimal prime of $J$.
By Theorem 5.2~(b) the ideal
$
I':=(\,I \,: \,b^{[d]} \,)
$
is binomial for suitably divisible $d$. Again, this quotient
is strictly larger than $I$ but $I'_{(J)}=I_{(J)}$,
so again we are done by Noetherian induction.
\Box

\vskip .2cm

As a corollary we deduce that
the minimal primary components of a binomial ideal
are all binomial.  Following
Eisenbud-Hunecke-Vasconcelos [1992], we write
$Hull(I)$ for the intersection of the minimal primary
components of an ideal $I$.
Note that $Hull(I) = I_{(\sqrt{I})}$.

\proclaim {Corollary 6.5}.  If $I\subset S$ is a binomial
ideal and $P$ is a minimal prime of $I$,
then the $P$-primary component of $I$ is binomial.  If
$I$ is a cellular binomial ideal, then $Hull(I)$ is
also binomial.

\noindent {\sl Proof}. By Theorem 6.2, we may assume
that $I$ is cellular for the
first statement, as well.  For the first
statement, take $J=P$ in
Theorem 6.4.  For the second statement, take $J = \sqrt{I}$
in Theorem 6.4.
\Box

\vskip .1cm

\proclaim Problem 6.6. Is $Hull(I)$ is binomial for
every (not necessarily cellular) binomial ideal~$I$~?

\vskip .6cm

\beginsection 7. Primary decomposition into binomial ideals

We are now in a position to do binomial primary
decomposition.

\proclaim Theorem 7.1. Let $I $ be a binomial ideal in
$ S= k[x_1,\dots,x_n]$.
If $k$ is an algebraically closed field
then $I$ has a minimal primary decomposition
in terms of binomial ideals.

The situation turns
out to be quite different in characteristic 0 and in characteristic
$p>0$. Curiously, the characteristic 0 result is far more
difficult and subtle.  We next formulate a
more precise result that makes the difference clear:

If $I$ is a binomial ideal in $S = k[x_1,\ldots,x_n]$,
then we write
${\cal Z}_I\subseteq \{1,\dots,n\}$ for the set
of indices $i$ such that $x_i$ is a nonzerodivisor
modulo $I$.  We write $M(I) = (\{x_i\}_{i\notin {\cal Z}_I})$
for the ideal generated by the other variables.
If the characteristic of $k$ is $p>0$
and $q = p^e$ is a power of $p$, then we write
$I^{[q]}$ for the ideal generated by
the $q^{\rm th}$ powers of elements of $I$.

\proclaim Theorem 7.1'. Let $I $ be a binomial ideal in
$ k[x_1,\dots,x_n]$, where $k$ is algebraically
closed.
\item{(a)} If $k$ has characteristic $p>0$ then,
for sufficiently large powers $q=p^e$,
$$
I \quad = \bigcap_{P\in Ass(S/I)}
 Hull\biggl(I+P^{[q]}\biggr) \eqno (7.1)
$$
is a minimal primary decomposition into binomial ideals.
\item{(b)} If  $k$ has characteristic 0, and $e$ is
a sufficiently large integer, then
$$
I \quad = \,\bigcap_{P\in Ass(S/I)}
Hull\biggl(I+M(P)^e+(P\cap k[{\cal Z}_P]) \biggr) \eqno (7.2)
$$
is a minimal primary decomposition into binomial
ideals.

\vskip .1cm

\noindent {\sl Remark: }  Formula (7.2) fails in positive
characteristic. For example, if
${\cal Z}_P= \{1,\dots,n\}$ for all  $P \in Ass(S/I)$
(the Laurent case),
then (7.2) states that $I$ is the
intersection of its associated primes, or, equivalently,
$I$ is radical. This is true only in
characteristic $0$.

\vskip .2cm

The reason why the positive characteristic case
is simpler is that in
positive characteristic we can make use of a variant of
the following result, which is part of the folklore of
primary decomposition:

\proclaim Lemma 7.2.
If $I$ is an ideal in a Noetherian ring $R$
and $e$ is a sufficiently large integer,
then $I$ has a minimal
primary decomposition of the form
$$ I \,\,\,\, = \bigcap_{P\in Ass(R/I)} Hull(I+P^e).
\eqno (7.3) $$

\noindent {\sl Proof: }
Let  $I = \cap Q_j$ be any minimal primary decomposition
of $I$, with $Q_j$ primary to a prime ideal $P_j$.
 The ideal  $P_j$ is the radical of $Q_j$, so
for large
$e$ we have
$P_j^e\subseteq Q_j$ for each $j$.
Thus
$Hull(I+P_j^e)$ is a $P_j$-primary
ideal containing $I$ and contained in $Q_j$.
It follows that $\,I = \cap_j Hull(I+P_j^e) \,$
is a primary decomposition of $I$.  It is minimal because
in the intersection (7.3)
there is only one ideal primary to each $P_j$.
\Box

If we start with a binomial ideal $I$
in a polynomial ring $S$, then
formula (7.3) does not give a binomial
primary decomposition because $P^e$ is almost never
binomial.  On the other hand,
the operation ``$Hull$'' preserves binomiality,
at least in the cellular case, by Corollary 6.5.
It would therefore
suffice to reduce to the cellular case and
replace $I+P^e$ in (7.3) by any binomial ideal
$J_P$ satisfying: \ \ \
(i) \ $I \subseteq J_P$, \ \ \
(ii) \  $\sqrt{J_P}=P$, \ \ \
(iii) \ $J_P \subseteq I + P^e$.

In characteristic $p>0$ it is easy to find
ideals satisfying (i),(ii), (iii), and hence to
prove Theorem 7.1.
But in characteristic 0, there is generally no binomial
ideal $J_P$ satisfying condition (i),(ii),(iii) for $e \geq 2 $.
(As an example take $I = (xy-x, x^2)$ and $P= (y-1,x)$.)
This accounts  for the difficulty of the characteristic 0
part of the proof given below.

\vskip .2cm

\noindent {\sl Proof of Theorem 7.1: }
(a)  \ Suppose that $P \in Ass (S/I)$.
By Theorem 6.1, $P$ is binomial.
The ``Beginner's Binomial Theorem''
 $\,(x+y)^q = x^q + y^q\,$ shows that $P^{[q]}$ is
binomial. For large  $e$ we have
$$\eqalign{
I \quad &\subseteq \,\bigcap_{P\in Ass(S/I)}
 Hull\biggl(I+P^{[q]}\biggr) \cr
&\subseteq \,\bigcap_{P\in Ass(S/I)}
 Hull\biggl(I+P^{q}\biggr) \quad = \quad I \cr}$$
by Lemma 7.2.
By Corollary 6.5 the terms
$Hull(I+P^{[q]})$
are binomial, and we are done.

\vskip .1cm

(b) If $P\in Ass(S/I)$ then
$P=M(P)+(P\cap k[{\cal Z}_P])$ is the only
minimal prime of $I+M(P)^e+(P\cap k[{\cal Z}_P])$,
so the hull of this ideal is primary.
Thus it suffices to prove equation (7.2).
 Let $I'$ denote the ideal
on the right hand side of equation
 (7.2).  Clearly, $I\subseteq I'$, we must prove
that $I'\subseteq I$.

Theorem 6.2
shows that $I$ is the intersection of cellular binomial
ideals $I_{\cal Z}^{(d)}$.  If we replace $I$ by
 $I_{\cal Z}^{(d)}$ then the terms
corresponding to associated primes of
 $I_{\cal Z}^{(d)}$ on the right hand  side
of (7.2)  do not change.  Thus it suffices to prove
$I'\subseteq I$  under the assumption that
$I$ is cellular.  In this case
$M(P)^e \subseteq I$ for $e \gg 0$, so we must show that
$$ I \quad \supseteq \quad  I' \,\,\, := \bigcap_{P\in
Ass(S/I)}  Hull\biggl(I+(P\cap k[{\cal Z}_P]) \biggr)
\eqno (7.2') $$
It suffices to prove this containment locally
at each prime $P$ in $Ass(S/I)$.
By Corollary 2.4 and Theorem 6.1, each such prime
can be written as  $P= I_+(\sigma)+M(P)$
for some partial character $\sigma$ on the lattice of
monomials in  $k[{\cal Z}_P]$. Set $K :=I_+(\sigma)S
= P \cap k[Z_P]$.

We shall do induction on the codimension of $P$ modulo
$I$.  We may therefore assume that $I'_{P'}=I_{P'}$
for all associated primes
$P'$ of $I$ properly contained in $P$.
It follows that $I'_P/I_P$ has support $P_P$, and thus has
finite length over the local ring $S_P$.
Equivalently, $I'_P \subseteq (I_P:P_P^\infty)$.
{}From the  definition of $I'$ we see that $I'_P \subseteq
Hull(I+K)_P$. Since  $I+K$ contains a power of $P$, we have
$Hull(I+K)_P=(I+K)_P$. Thus
$I'_P \subseteq (I_P:P_P^\infty)\cap (I+K)_P$,
and it suffices to show that
$(I_P:P_P^\infty)\cap (I+K)_P\subseteq I_P$.  Equivalently,
it suffices to show
$$ (I:P^\infty)\cap K \quad \subseteq \quad I_P. $$

By Theorem  6.4 the ideal $I_{(P)} = I_P\cap S$
is binomial. Replacing $I$ by $I_{(P)}$
does not change $I_P$, and at worst makes the ideal
$ (I:P^\infty) $  larger.  Thus we may suppose that $I =
I_{(P)}$ from the outset -- that is, we may suppose that
$P$ is the unique
maximal associated
prime of $I$. In this situation we shall finish the
proof by showing that
$$ (I:P^\infty)\, \cap \, K \quad \subseteq \quad I.
\eqno (7.4) $$

We fix a sufficiently divisible integer $d$.
Let $\sigma_d$ be the
restriction of $\sigma$ to the sublattice
$dL_\sigma \subset L_\sigma$. Since $P$ is the unique
maximal associated prime of $I$,
the saturations of $\sigma_d$ other than $\sigma$
correspond to primes $P_i$ that contain nonzerodivisors
modulo $I$. Since the characteristic of $k$ is 0, Corollary
2.2 shows that
$I_+(\sigma_d) = P\cap\bigcap_i P_i$.  Thus
$$
(I:I_+(\sigma)) \, \subseteq\,
(I:I_+(\sigma_d)) \,=\,
(I:P\cap\bigcap_i P_i) \,\subseteq\,
(I: I_+(\sigma)\prod_i P_i)\, =\,
(I:I_+(\sigma)).
$$
It follows that
$(I:I_+(\sigma))= (I:I_+(\sigma_d))$.

Corollary 5.6 shows that   $(I:I_+(\sigma_d))=
(I:I_+(\sigma_d)^\infty)$. We have $K \subset P$
and $P^N \subset I + K $ for $N \gg 0$.
Putting everything together we get
$$\eqalign{
(I:P^\infty) \,\,&\supseteq \,\,(I:K ^\infty)
\,\, \supseteq \,\,
\,\, (I:K) \,\, = \,\, (I:I_+(\sigma)) \,\,\cr
&= \,\, (I:I_+(\sigma_d))\,\,=
\,\,(I:I_+(\sigma_d)^\infty)\,\,
\supseteq \,\,
(I:I_+(\sigma)^\infty) \,\, = \,\,(I:P^\infty). }$$
In particular $(I:P^\infty)=(I:K) $, and hence
(7.4) is equivalent to $(I:K) \cap K \,\subseteq I$.

By Theorem 5.2 (b) applied to the generators
$b^{[d]}$ of $I_+(\sigma_d)$, there exists a monomial
ideal $J$ such that $(I:K) =I+J$. The variables
$x_i$ for $i\in {\cal Z}$ are nonzerodivisors modulo
$I$ hence also modulo $(I:K)$. Therefore the minimal
generators of $J$ may be taken to be
 monomials in the other variables
$\{x_i\}_{i\notin {\cal Z}}$ alone.
Suppose that $f \in K\cap (I:K). $
We may write $f = g+j$ where $g \in I$
and $j\in J$.  We shall show that $f \in I$ by
an induction on the number of terms of $j$, the
case $j= 0$ being trivial.

Consider any reduced  Gr\"obner basis ${\cal G}$ of $K$. By Corollary 2.3,
${\cal G}$ consists of binomials $x^{u_+} - \sigma(u) x^{u_-}$
where $u \in L_\sigma \subset {\bf Z}^{\cal Z} \subset
{\bf Z}^n$. Reduction of a monomial $x^\alpha$ modulo ${\cal G}$
results in a nonzero term whose exponent vector is congruent to
$\alpha $ modulo the lattice $L_\sigma$.

Let $j_1$ be a term of $j$, and $m$ a generator
of $J$ that divides $j_1$.
Since  $f$ reduces to zero mod ${\cal G}$, we see that there is another
term $j_2$ of $f$ whose exponent vector is congruent to
that of  $j_1$ modulo $L_\sigma$. In particular,
since none of the variables dividing
$m$ is in $k[{\cal Z}]$, the term $j_2$ is also divisible
by $m$, and for some constant $c\in k^*$ the binomial
$(j_1/m)-c(j_2/m)$ is in $K$.  As $m\in (I:K)$,
we have $j_1-cj_2 \in I$.  Subtracting  $j_1-c j_2$
from $f$, we get a new element $f' \in K\cap (I:K) $
with a decomposition $f' = g'+j'$,
where $g' \in I, j' \in J$, and
$j'$ has fewer terms than $j$.
By induction, $f'\in I$ and therefore $f \in I$.
\Box

\vskip .2cm

In spite of Theorem 7.1 there are still many open questions
about the decomposition of binomial ideals.  For example:

\proclaim Problem 7.3. \rm
Does every binomial ideal have an irreducible
decomposition into binomial ideals ?
Find a combinatorial characterization of
irreducible binomial ideals.

\vskip .2cm

If we are given any ideal $I $ in $S$, then a prime ideal
$P$ is associated to $I$ if and only if there exists
$f \in S$ such that $(I:f) = P$. Such polynomial $f$
might be called a {\it witness} for the prime $P$.
In the case where $I$ is binomial and hence $P$
is binomial, one may ask whether there exists
a binomial witness. The answer to this question
is easily seen to be ``no'': take
$P = (x-1)$ and $ I = (x^d-1)$, where every witness,
like  $1+x+\dots + x^{d-1}$, has at least $d$ terms; or take
$P = ( x_1,x_2,...,x_n)$, the ideal of the origin,
and  $I = ( \{x_i^2 - x_i\}_{i =1,\dots,n})$,
the ideal of the vertices of the cube,
where it is easy to show that
any witness, like $\prod (x_i-1)$, has at least $2^n$ terms.
However, the following ``Witness Theorem'' provides a
monomial witness in a restricted sense:

\proclaim Theorem 7.4. \ Let
$I$ be a cellular binomial ideal in
$S = k[x_1,\ldots,x_n]$, and let ${\cal Z} = {\cal Z}(I)$.
If $\,P \,=\, I_+ (\sigma) + M ({\cal Z}) \,$
is an associated prime of $I$, then
there exists a monomial $m$ in
the variables
$\{x_i\}_{i\notin {\cal Z}}$ and
a partial character $\tau$ on ${\bf Z}^{\cal Z}$ such that
$\sigma$ is a saturation of $\tau$ and
$$\,( I \, : \, m ) \,\cap \,k[ {\cal Z} ]
\quad = \quad  I_+ (\tau) .$$

\noindent {\sl Proof: }
The proof is by Noetherian induction. First,
if $I$ contains all the variables $\{x_i\}_{i\notin {\cal Z}}$,
then we are in the Laurent case: $\, I \,=\, I_+(\tau)
+ M({\cal Z})\,$ for some $\tau$, by Corollary 2.3.
In this case the assertion holds with $m=1$.
 Otherwise there exists a variable, say
$x_1$ after relabeling, such that both the cellular
ideals $(I : x_1)$ and
$$ \,I' \quad :=\quad \bigl(\,( I + (x_1)  ) :
(\prod_{i\in {\cal Z}}x_i)^\infty \,\bigr) \,
$$
strictly contain $I$.

By  Noetherian induction we may assume
that Theorem 7.4 holds for $(I:x_1)$ and $I'$.
As in the proof of Theorem 6.1, every associated
prime $P$ of $I$ is associated to $(I: x_1)$ or to $I'$.
If $P$ is associated to $(I:x_1)$,
then we have a presentation
$$\,( \,( I : x_1)  \, : \, m' ) \,\cap \,
 k[ {\cal Z} ] \quad = \quad  I_+ (\tau) $$
for some monomial $m'$.  Taking $m = x_1 m'$,
the claim follows.

We may therefore assume that $P$ is associated to $I'$.
By the Noetherian induction again, there exists a monomial
$m$ and a partial character $\tau$ with saturation
$\sigma$ such that
$$
\biggl( \bigl(( I + (x_1)  ) :
 (\prod_{i\in {\cal Z}}x_i)^\infty \bigr) \,:\,
 m  \biggr) \
 \,\,\cap \,\,k[{\cal Z}] \quad = \quad  I_+ (\tau) .
\eqno (7.5)
$$
We claim that this ideal equals $( I  :  m ) \,\cap \,k[{\cal Z}]$.
Certainly $( I  :  m ) \,\cap \,k[{\cal Z}]$ is contained in (7.5).
Note also that (7.5) is a proper ideal.

Let $f$ be any polynomial in (7.5).
Suppose that $mf$ has a term in  $I + (x_1)$.
Since the terms in $f$ are all in
$k[{\cal Z}]$, we would have $\,m \in ((I + (x_1)):
(\prod_{i \in {\cal Z}} x_i)^\infty)$,
and the ideal in (7.5) would not be proper.  Therefore no term of $mf$
is in $I + (x_1)$. Using Proposition 1.10 we conclude that
$mf \in I$, as required.
\Box

\vskip .2cm

Using Theorem 7.4, we get the following alternative
decomposition of a binomial ideal.
We conjecture that Corollary 7.5 holds
in finite characteristic as well.

\proclaim Corollary 7.5. Let $k$ be a field of
characteristic $0$, let $I$ be a cellular binomial
ideal in $S= k[x_1,\ldots,x_n]$,
and ${\cal Z} = {\cal Z}(I)$. Then
$I$ has the following presentation as a finite intersection
of unmixed binomial ideals:
$$
I \quad = \! \bigcap_{m\ {\rm a\ monomial\ in}\
\{x_i\}_{i\notin {\cal Z}}} \!\!\!\!
Hull \biggl(I+\bigl((I:m)\cap k[{\cal Z}]\bigr)\biggr).
\eqno (7.6)
$$

\noindent {\sl Proof:}
The intersection given in (7.6) clearly contains $I$.
On the other hand, if $P = I_+(\sigma) + M({\cal Z})$
is an associated prime of $I$ then by
Theorem 7.4 there is a monomial $m$ in the variables
$\{x_i\}_{i\notin {\cal Z}}$ such that
$(I:m)\cap k[{\cal Z}] = I_+(\tau)$,
and $\sigma$ is a saturation of
$\tau$.  Thus
$$
Hull (I+((I:m)\cap k[{\cal Z}])) \, = \,
Hull (I\,+\, I_+(\tau) ) \, \subseteq \,
Hull (I\,+\, I_+(\sigma) ) \, =\,
Hull (I+(P\cap k[{\cal Z}]));
$$
hence the intersection in formula (7.6) is contained in the
intersection in formula (7.2).
\Box

\vskip .2cm

The first step in the computation of a
primary decomposition of a binomial ideal $I$
is to find a cellular decomposition as in (6.3).
In certain cases the cellular decomposition is already a
primary decomposition. We next show that this event happens
when the algebraic set defined by $I$ is irreducible
and not  contained  in any coordinate hyperplane.

\proclaim Theorem 7.6.
Let $I \subset S $ be a binomial ideal.  If
 $\sqrt{I}$ is prime and does not contain
any of the variables, then for all
sequences $d$ of sufficiently large integers
the ideals $I_{\cal Z}^{(d)}$
are primary.  Thus the cellular
decomposition (6.3)
is a (possibly nonminimal) primary decomposition of $I$.

\vskip .1cm

\noindent {\sl Proof: }
Set $P =\sqrt{I}$, and let ${\cal Z}$ be any
subset of $\{1,\ldots,n\}$.
We define $I_{\cal Z}^{(d)}$
as in formula (6.3) and $P_{\cal Z}$ as in formula (4.2).
Clearly, $I_{\cal Z}^{(d)} \subseteq
P_{\cal Z} \subseteq
\sqrt{ I_{\cal Z}^{(d)} } $, so that
 $I_{\cal Z}^{(d)}$ is a proper ideal if and only if
$P_{\cal Z}$ is a proper ideal. In this case,
$P_{\cal Z} \cap k[{\cal Z}] =P\cap k[{\cal Z}]$, by Lemma 4.3,
and thus $ P_{\cal Z} = (P\cap k[{\cal Z}]) + M({\cal
Z})$. This shows that $P_{\cal Z}$
is prime, so  $P_{\cal Z} = \sqrt{I_{\cal Z}^{(d)}} $.
We conclude that every associated prime
of $I_{\cal Z}^{(d)}$ contains $P_{\cal Z}$.

Let $Q $ be any associated prime
of $I_{\cal Z}^{(d)}$. By Theorem 6.1 and Corollary 2.4, we can write
$Q = I_+(\sigma)S+M({\cal Z})$, where $\sigma$
is a partial character on ${\bf Z}^{\cal Z}$.
By Theorem  7.4, there exists a positive integer $e$
and a monomial $m\notin I_{\cal Z}^{(d)}$ such that
$I_+(\sigma_e)m \subset I_{\cal Z}^{(d)}$.
(Here $\sigma_e$ denotes the restriction of
$\sigma$ to the lattice $e L_\sigma \subseteq L_\sigma$.)

Let $f$ be any element of $I_+(\sigma_e)$.
Then $f m \in I_{\cal Z}^{(d)} $, so there exists
a monomial $m'$ in $k[{\cal Z}]$ such that
$fmm' \in I+(\{x_i^{d_i}\}_{i\notin {\cal Z}})$.
Since $mm' \notin I_{\cal Z}^{(d)}$ and
$f\in k[{\cal Z}]$, the terms of $fmm'$ are not in
$I+(\{x_i^{d_i}\}_{i\notin {\cal Z}})$.
It follows by Proposition 1.10 that
$fmm' \in I \subseteq P $.
Since the prime $P$ does not contain any monomials,
it follows that $f \in P$.  This shows that
$I_+(\sigma_e)$ is contained in $ P$.   Since $P\cap
k[{\cal Z}]$ is contained in $I_+(\sigma)$, it follows
that $I_+(\sigma) = P\cap k[{\cal Z}]$ and consequently $Q
= P_{\cal Z}$. We conclude that $P_{\cal Z}$ is the only
associated prime of $I_{\cal Z}$.
\Box

\vskip .2cm

\noindent {\bf Example 7.7.}
Theorem 7.6 does not hold in general
for binomial ideals $I$ whose
radical is prime but
contains a variable $x_i$.
For example,
 $\, I \, = \, \bigl(x_1^2, x_1(x_2-x_3) \bigr)
\, = \, (x_1) \, \cap \, (x_1^2, x_2-x_3)\,$
has radical $(x_1, x_2-x_3)$, but if
${\cal Z}=\{ 2,3 \}$
then $I = I_{\cal Z}^{(d)}$ is not primary.

\vskip .3cm

We next determine
which of the ideals $P_{\cal Z}$
arising in  the proof of Theorem 7.6 is proper (this
is somewhat weaker than saying that the
corresponding cell $(k^*)^{\cal Z}$ contains
an associated point of $I$).
This condition is phrased in terms of combinatorial
convexity. It is well-known in the theory of toric
varieties. Let $P = I_+(\sigma)$ be a binomial prime ideal
in $S $ such that $x_i \not\in P$ for all $i$. Let $d = dim(P)$.
Then ${\bf Z}^n/ L_\sigma $ is a free abelian
group of rank $d$, and
$V = ({\bf Z}^n/ L_\sigma ) \otimes_{\bf Z} {\bf R}$
is a $d$-dimensional real vector space.
Let $\bar e_i$ denote the image in $V$
of the $i$-th unit vector in ${\bf Z}^n$.
We consider the $d$-dimensional convex polyhedral cone
$$ {\cal C} \quad := \quad
\bigl\{\, \lambda_1 \bar{e_1} +
\lambda_2 \bar e_2  + \cdots
+ \lambda_n \bar e_n \,\,:\,\,
\lambda_1,\lambda_2,\ldots,\lambda_n \geq 0 \,\bigr\}. \eqno (7.7) $$
A subset ${\cal Z}$ of $\{1,\ldots,n\}$
is said to be a ${\it face}$ of $P$ if
$\,pos(\{\bar e_i\,:\,i \in {\cal Z}\}) \,$
is a face of ${\cal C}$.

\proclaim Proposition 7.8.
With notation as above,
the ideal $P_{\cal Z}$ is proper if and only if
${\cal Z}$ is a face of $P$.

\noindent {\sl Proof: }
Suppose ${\cal Z}$ is not a face. By elementary
convexity, this is equivalent to the following: the
generators of  ${\cal C}$ satisfy a linear dependency
of the form $\, \lambda_1 \bar e_{i_1} + \cdots +
    \lambda_s \bar e_{i_s} \, = \,
   \mu_1 \bar e_{j_1} + \cdots  + \mu_t \bar e_{j_t}$,
where $\lambda_1,\ldots,\lambda_s,
\mu_1,\ldots,\mu_t$ are positive integers,
$\{i_1,\ldots,i_s \} \subseteq {\cal Z}$, and
$\{j_1,\ldots,j_t \} \not\subseteq {\cal Z}$.
The ideal $P$ therefore contains some binomial
$\, x_{i_1}^{\lambda_1} \cdots    x_{i_s}^{\lambda_s} \,
 - \, c \cdot   x_{j_1}^{\mu_1} \cdots x_{j_t}^{\mu_t}$,
$\, c \in k^* $. This shows that a power of
$\, x_{i_1}^{\lambda_1} \cdots    x_{i_s}^{\lambda_s}
\,$ lies in  $\, P + M( {\cal Z}) $, and
consequently $P_{\cal Z}$ contains a unit.
 Conversely, let ${\cal Z}$ be a face.
Then there is no linear dependency as above, which means
that every binomial in $P$ lies in $k[{\cal Z}]$
or in $M({\cal Z})$. Therefore
$P_{\cal Z} = (P\cap k[{\cal Z}])+ M({\cal Z})$,
and this is  clearly a proper ideal.
\Box

\vskip .3cm

Proposition 7.8 can be rephrased as follows.
If an ideal $I$ satisfies the
hypothesis of Theorem 7.6, then  its
associated points are in natural bijection
with a subset of the faces of $\sqrt{I}$.
We close this section by describing a
class of binomial ideals with these properties.

\vskip .2cm

\noindent {\bf Example 7.9.} \  {\sl (Circuit Ideals)}
\ \ \
Let $\rho$ be a saturated partial character on
$\,{\bf Z}^n $.  If $v \in {\bf Z}^n $, then the
{\it support} of $v$ is the set of basis elements of
${\bf Z}^n $ that appear with nonzero coefficient in
the expression of $v$.
A primitive non-zero element $v$ of the lattice $L_\sigma$ is
said to be a {\it circuit} if the support of $v$ is
minimal with respect to inclusion.
The {\it circuit ideal}  $C(\rho)$ is the ideal
generated by the binomials
$x^{\alpha_+} - \rho(\alpha) \,x^{\alpha_-}$,
where $\alpha$ runs over all circuits of $L_\sigma$.
Clearly, $C(\rho)$ is contained in the prime ideal
$I_+(\rho)$.
For certain special lattices $L_\rho$ arising in combinatorics
we have $C(\rho) = I_+(\rho)$; for instance, this
is the case for lattices presented by totally
unimodular matrices \ (see
\S 4 of (Sturmfels [1992])). In general we have
only:

\proclaim {Proposition 7.10}.  With notation
as above,
$$\,\sqrt{C(\rho)} \,= \,I_+(\rho).$$

In particular, we see that Proposition
7.8 applies to circuit ideals.
For the proof,
we need to know that $L$ is generated by circuits,
which is a special case of the following:

\proclaim Lemma 7.11. Let $R$ be an integral domain.
If $\phi: R^n \to R^d$ is an epimorphism, then the
kernel of $\phi$ is the image of the map
$\,\psi \,: \,\wedge_{d+1} R^n\, \rightarrow\,
 R^n, \,\, \xi \,\mapsto \, \xi \,\lfloor\,
\wedge_d \phi $.  The circuits in the kernel of $\phi$
are,  up to multiplication by elements
of the quotient field, precisely the
nonzero images of the standard
basis vectors of $\wedge_{d+1} R^n$.  These
images are the
relations given by Cramer's rule,
$$ \psi \,( e_{i_0} \wedge e_{i_1} \wedge
\cdots \wedge e_{i_d} ) \,\,\,\, = \,\,\,\,
\sum_{j=0}^d \,(-1)^j \cdot det\,(\phi_{i_0,\ldots,i_{j-1},
i_{j+1},\ldots,i_d}) \cdot e_{i_j} . $$

\noindent {\sl Proof: }
To prove the first statement,
let $U$ be a $d \times n$-integer
matrix such that $\phi \cdot U $ is the $d \times d$-identity
matrix. If $v \in ker(\phi)$, then an elementary computation
in multilinear algebra gives:
$$ \psi(\wedge_d U \wedge v ) \quad = \quad
(\wedge_d U \wedge v ) \,\lfloor\, \wedge_d \phi
\quad = \quad  ((\wedge_d \phi) \cdot (\wedge U)) \cdot v
\quad = \quad v .$$

Call the relations $\psi \,( e_{i_0} \wedge e_{i_1} \wedge
\cdots \wedge e_{i_d} )$ {\sl Cramer relations}.
If a Cramer relation is nonzero, then it is a relation
among $d+1$ elements of $R^d$
that generate a submodule of rank $d$ in $R^d$.  Any relation
among these $d+1$ images must be a multiple of the Cramer
relation by an element of the quotient field.  In
particular, the Cramer relation is a circuit in $ker(\phi)$.

To show conversely that every circuit
in $ker(\phi)$ is,
up to multiplication by an element of the
quotient field, a Cramer relation, it
now suffices to prove that every circuit has
support contained in a set of $d+1$ vectors
whose images generate a module of rank $d$.
For every circuit $c$ there is
a number $r$ such that $c$ is a relation among
$r+1$ vectors spanning a submodule
of rank exactly $r$ (if the rank were lower,
then there would be two independent relations,
and thus a relation involving a subset of the
terms).  Since the images of the basis vectors
of $R^n$ span a submodule of rank $d\geq r$,
we can find $d-r$ such vectors whose images,
together with the images of the vectors in the
support of $c$, span a submodule of rank $d$,
and we are done.
\hfill \Box

\vskip .2cm

\noindent {\sl Remarks:}
The statements about circuits
are false if $R$ is not an
integral domain:  if $x$ is a zerodivisor, then the only
circuits in the kernel of the
map $(1,x): R^2\to R$ are the column
vector with entries $0,y$, where $xy=0$,
so the relation defined
by Cramer's rule is not a circuit, and the
circuits do not generate all the relations.
However the Cramer relations still do
generate, and this fact
has been extended by Buchsbaum and Rim [1964]
to a natural free resolution.

\vskip .3cm

\noindent {\sl Proof of Proposition 7.10.}
It is easy to see that
every element of $L_\rho$ is a positive rational
linear combination of circuits. Therefore the convexity
argument in the proof of Proposition 7.8
applies to circuit ideals as well, and $C(\rho)_{\cal Z}$
is a proper ideal if and only if ${\cal Z}$ is a face of
$I_+(\rho)$. Now, suppose that $\,{\cal Z} \subset
\{1,\ldots,n\}\,$ is a face. Let $\,\rho|_{\cal Z} \,$
denote the restriction of $\rho$ to the sublattice
$\,L_\rho \, \cap \,{\bf Z}^{\cal Z}$. By Lemma 7.11
applied to this sublattice, we have $\,I_+(\rho|_{\cal Z})
\, = \,(C(\rho|_{\cal Z}) :
 (\prod_{i \in {\cal Z}} x_i)^\infty) $.
Clearly, the circuits of  $\,L_\rho \, \cap \,
{\bf Z}^{\cal Z}$ are just the circuits of $L_\rho$
that have support in ${\cal Z}$. Hence
$\,C(\rho|_{\cal Z})\,\subseteq \,
C(\rho) \cap k[{\cal Z}]  \,$
and we conclude
$$  \eqalign{  C(\rho)_{\cal Z} \,\,\, & = \,\,\,
\biggl( \bigl( C(\rho) + M({\cal Z}) \bigr) \,:\,
(\prod_{i \in {\cal Z}}\! x_i)^\infty \biggl)
\,\,\, = \,\,\,
\biggl( \bigl( ( C(\rho) \cap k[{\cal Z}])
+ M({\cal Z})\bigr) \,:\, (\prod_{i \in {\cal Z}}\!
x_i)^\infty \biggl)  \cr & \supseteq \,\,\,
\bigl( \,C(\rho|_{\cal Z}) :
(\prod_{i \in {\cal Z}}\! x_i)^\infty \,\bigl)
\,\,+ \,\,M({\cal Z})
 \,\,\,\, =\,\, \,\, I_+ (\rho|_{\cal Z}) \,+ \,M({\cal
Z})  \, \,\, = \,\,\,I_+ (\rho)_{\cal Z} .\cr} $$
Since the reverse inclusion is obvious, we have
$C(\rho)_{\cal Z} = I_+(\rho)_{\cal Z}$.
Our claim follows by taking the intersection
over all faces ${\cal Z}$ of $ I_+ (\rho)$.
\Box

\vskip .2cm

\noindent {\bf Problem:}
It remains  an interesting combinatorial
problem to characterize the embedded primary
components of the circuit ideal $C(\rho)$. In particular,
which faces of (the polyhedral cone associated with
the prime) $\, I_+(\rho)\,$
support an associated prime of
$C(\rho)$~? \ An answer to this question might be
valuable for the applications of binomial ideals
to integer programming and statistics
mentioned in the introduction.

\vskip .6cm

\beginsection 8. Algorithms

In this final section we present algorithms
for computing various aspects of the
primary decomposition of a binomial ideal.
In each case we outline only the basic steps, and
we disregard questions of efficiency. It remains
an  interesting problem to find
best possible procedures. Our Algorithms 8.1 -- 8.7
differ greatly from the known algorithms for general
polynomial ideals, given for example by
Gianni-Trager-Zacharias [1988]  and
Eisenbud-Hunecke-Vasconcelos [1992]. The older algorithms
immediately leave the category of binomial ideals
(in the sense that they either make changes of coordinates
or use syzygy computations and Jacobian ideals).
The algorithms presented below work almost entirely
with binomials and thus maintain maximal sparseness.
This is an important advantage because sparseness is a significant
factor in the effectiveness of computations.

\vskip .3cm

\noindent {\bf Algorithm 8.1: Radical. }

\noindent {\sl Input: \ } A binomial ideal $I$ in
$S = k[x_1,\ldots,x_n]$.

\noindent {\sl Output: \ } A finite set of
binomials generating the radical $\sqrt{I}$ of $I$.

\item{1.} If $I=(0)$ output $\{0\}$.
If $I= S$ output $\{1 \}$.
\item{2.} Otherwise compute
$\, J \, = \,(I:(x_1 \cdots x_n)^\infty)$, for instance
by introducing a new variable $t$ and eliminating $t$ from
$\,I \,+\, (\,t \,x_1 \cdots x_n - 1\,)$.
\item{3.} If $char (k) = p>0$ compute the radical of
$J$ by computing the $p$-saturation of its
associated lattice as in Corollary 2.2.
Set $J := \sqrt{J}$.
\item{4.} For $i = 1,\ldots,n $ do
\itemitem{4.1} Replace $x_i$ by $0$ in all
generators of $I$. \hfill \break
Let $J_i$ be the resulting
ideal in  $k[x_1,\ldots,x_{i-1},x_{i+1},\ldots,x_n]$.
\itemitem{4.2}
Compute $\,\sqrt{J_i} \,$ by recursively calling Algorithm 8.1.
\item{4.}
Compute and output a reduced Gr\"obner basis
for the intersection
$$ J \,\cap \,
\, \bigl(\sqrt{J_1} S  +  (x_1) \bigr) \,\cap \,
\, \bigl(\sqrt{J_2} S  +  (x_2) \bigr)
\,\cap \,\,\cdots \,\,\cap
\, \bigl(\sqrt{J_n} S  +  (x_n) \bigr) . $$

\noindent {\sl Comments: }
The correctness follows
from the results in Sections 2 and 3,
in particular Theorem 3.1 and formula (3.1).
If the characteristic of $k$ is $0$ then Step 3 is
unnecessary: in this case $J$ is already
radical by Corollary 2.2. As it stands
Algorithm 8.1 requires $n ! $ recursive calls.
The following algorithm accomplishes the same
task in $2^n$ iterations.

In what follows we use the  abbreviation
$M := M({\cal Z}) = (\{ x_i \}_{i \not\in {\cal Z}})$.

\vskip .5cm

\noindent {\bf Algorithm 8.2: Minimal primes.}

\noindent {\sl Input: \ } A binomial ideal $I$
in $S = k[x_1,\ldots,x_n]$.

\noindent {\sl Output:} Binomial prime
ideals $P_1,\ldots,P_s$ whose intersection
is irredundant and equals~$\sqrt{I}$.

\item{$\bullet$} For each subset ${\cal Z} $ of
$\{1,\ldots,n\}$ do

\itemitem{1.} Decide whether the ideal
$$ I_{\cal Z} \quad := \quad \bigl(
( \,I \, + \, M({\cal Z}) \,) \,: \,
( \prod_{i \in {\cal Z}} x_i )^\infty \, \bigr)$$
is proper. If not, stop here. Otherwise continue.
\itemitem{2.} Determine the unique partial character
$\rho $ on ${\bf Z}^{\cal Z}$ such that
$\,I_{\cal Z} \, = \,I_+ (\rho) \,+ \, M({\cal Z})$.
\itemitem{3.} If $char(k) = p$ then
replace $\rho$ by the unique extension $\rho'$
of $\rho$ to  the $p$-saturation of $L_\rho$ as in
Corollary 2.2.
\itemitem{4.} Compute the saturations
$\rho_1,\ldots,\rho_g$ of $\rho$, and save the $g$ primes
$\,I_+(\rho_i) + M({\cal Z})$.
\item{$\bullet$} Among all prime ideals computed remove the
redundant ones and output the others.

\vskip .3cm

\noindent {\sl Comments: } The correctness
of Algorithm 8.2 follows from Theorem 6.1 and
the results in Section 4. In the worst
case each of the $2^n$ subsets ${\cal Z}$
will contribute a minimal prime: this happens for
$ (\{x_i^2 - x_i , \,i=1,\ldots,n\})$ as in Example 4.5.
On the other hand, for many binomial ideals
we can avoid having to inspect all $2^n$ cells.
One natural shortcut arises if (in the course of the
algorithm) we find that
$I_{{\cal Z}_1} = S$ and $I_{{\cal Z}_2} \not= S$
for ${\cal Z}_1 \subset {\cal Z}_2 $. Then
we may ignore all subsets ${\cal Z}$
with ${\cal Z} \,\cap \,{\cal Z}_1 \,=\,{\cal Z}_2$:
for such ${\cal Z}$ the ideal $I_{\cal Z}$ cannot be
proper by Theorem 4.1 (ii). Also Proposition 7.8 allows
some savings: if $I$ is a cellular radical ideal, then
one may precompute the faces of the cone ${\cal C}$
in formula (7.7) using some convex hull algorithm.
The same techniques can be used to speed
up the next algorithm. The correctness of Algorithm 8.3
is essentially the content of Theorem 6.2.

\vskip .2cm

\noindent {\bf Algorithm 8.3: Cellular decomposition. }

 \noindent {\sl Input: \ } A binomial ideal $I$ in
$S = k[x_1,\ldots,x_n]$.

\noindent {\sl Output: }
Cellular ideals $J_{\cal Z}$, indexed by
a set of subsets of $\{1,\ldots,n\}$, such that
$\bigcap_{\cal Z} J_{\cal Z} = I$.

\item{1.} Fix a vector $d = (d_1,\ldots,d_n)$ of
sufficiently large integers \ (see Problem 6.3).

\item{2.} For each ${\cal Z} \subseteq \{1,\ldots,n\}$,
let $J_{\cal Z} := I^{(d)}_{\cal Z}$ as defined
in formula (6.3).

\item{3.} Output those proper ideals $J_{\cal Z}$
which are minimal with respect to inclusion.

\vskip .4cm

\noindent
In the remaining four algorithms we shall  restrict ourselves to
ideals which are cellular (i.e., $I = I^{(d)}_{\cal Z}$ for $d_i \gg 0$).
For general ideals this requires to run Algorithm 8.3 beforehand.

\vskip .2cm

\noindent {\bf Algorithm 8.4: Test for primary ideals. }

\noindent {\sl Input:} A subset ${\cal Z}\subset
\{1,\ldots,n\}$ and a binomial ideal $I$ which is
cellular with respect to ${\cal Z}$.

\noindent {\sl Output: \ } The radical of $I$,
and the decision (``YES'', ``NO'') whether $I$ is
primary. In the negative case the algorithm generates two
distinct associated prime ideals of $I$.

\item{1.} Compute the unique partial character
$\,\sigma \,$ on ${\bf Z}^{\cal Z}$ such that
$\,I \,\cap \, k[{\cal Z}] \, = \, I_+ ( \sigma) $.
If the characteristic is $p>0$, replace $\sigma$ by
its $p$-saturation.
\hfill \break
Output: ``The radical of $I$ equals
$I_+(\sigma) + M $''.

\item{2.} If $\sigma$ is not saturated, then output ``NO,
the radical of $I$ is not prime'', choose two
distinct saturations $\sigma_i$
and $\sigma_j$ of $\sigma$,
output the two associated primes
$I_+(\sigma_i) + M $ and
$I_+(\sigma_j) + M $, and STOP.

\item{3.} Compute a Gr\"obner basis of $I$, and let
${\cal T}$ be the set of {\it maximally standard}
monomials in the variables $\{x_i\}_{i\notin {\cal Z}}$.
In other words, ${\cal T}$ is equal to the
set of monomials in $\,(in(I) : M) \setminus ( in(I) +
(x_j ,\, j \in {\cal Z}) )$.

\item{4.} If   $\,(I : m ) \,\cap \, k [{\cal Z}] \, \subseteq \,
 I_+ (\sigma)$ for all $m \in {\cal T}$, then output
``YES, $I$ is primary''.

\item{5.}
Otherwise, choose $ m \in {\cal M}$  such that
$\,(I : m) \,\cap \, k [{\cal Z}] \, = \,I_+ ( \rho) \,
\not \subseteq \, I_+ ( \sigma) \ .$
Let $\rho'$ be
any saturation of $\rho$. \ \
Output: ``NO, $I$ is not primary. The primes
$I_+(\sigma) + M $ and
$I_+(\rho') + M $ are both associated to $I$.''

\vskip .2cm

\noindent {\sl Comments: } \  In light of Theorem 7.4, every
associated prime of $I$ is associated to \break
$((I:m)\cap k[{\cal Z}])+M$ for some monomial $m$
in the variables $\{x_i\}_{i\not\in {\cal Z}}$.
The maximal proper ideals of the form
$(I:m)$ are gotten from monomials $m$ in the
finite set ${\cal T}$ constructed in step 3.
In step 5,  the ideal $I_+(\rho)$ properly
contains the prime ideal $I_+(\sigma)$.
Therefore $I_+(\sigma)$ is properly contained
in any associated prime  $I_+(\rho')$ of $I_+(\rho)$.

\vskip .4cm

\noindent {\bf Algorithm 8.5: Associated primes}

\noindent {\sl Input:} A subset ${\cal Z}\subset
\{1,\ldots,n\}$ and a binomial ideal $I$ which is
cellular with respect to~${\cal Z}$.

\noindent  {\sl Output: \ } The  list of associated
primes $P_1,\ldots,P_s$ of $I$.

\item{1.} Compute a Gr\"obner basis of $I$.

\item{2.} Let ${\cal U}$ be the set of standard
monomials in the variables $\{x_i\}_{i\notin {\cal Z}}$.

\item{3.} For each $m \in {\cal U}$ do
\itemitem{3.1.} Compute the partial character $\tau$ that
satisfies $\,I_+(\tau) \,= \,(I: m) \,\cap \,k[{\cal Z}]$.
\itemitem{3.2.} For each saturation $\tau'$ of $\tau$
output the prime ideal $\, I_+ (\tau') + M $.

\vskip .2cm

\noindent {\sl Comments: }  The standard monomials
in step 2 are those not contained in the initial
ideal of $I$. The primes $I_+(\tau')+M$ are
associated to $I_+(\tau)+M$.  It follows that
$I_+(\tau')+M$ is associated to $I$.  Theorem 7.4 shows
that every associated prime of $I$ occurs in this way.
The set ${\cal U}$ is finite because a
power of $M = (x_i, \, i \not\in {\cal Z})$ lies $I$.
Note that the set ${\cal T}$ in step 3 of Algorithm 8.4
consists precisely of the maximal (with respect to
divisibility) monomials in ${\cal U}$.

\vskip .4cm

\noindent {\bf Algorithm 8.6: Minimal primary component. }

 \noindent {\sl Input: \ } A cellular binomial ideal $I$
whose radical $\sqrt{I}$ is prime.

\noindent {\sl Output: \ } A set of binomial generators
for  the primary ideal $Hull(I) \, = \,I_{(\sqrt{I})}$.

\item{0.} Set $J = \sqrt{I}$ and let $\sigma$ be the
saturated partial character such that $J = I_+(\sigma)+M$.

\item{1.} Call Algorithm 8.4 to determine whether
$I$ is primary. If yes, output $I$ and STOP. \break
If no, we get another associated prime
$P = I_+(\rho) + M$ properly contained in $J$.

\item{2.} We shall now follow the proof of
Theorem 6.4 verbatim. \hfill \break
First, compute the lattice $L$ in formula (6.5).

\item{3.} If $L$ has finite index in $L_\rho$ then
proceed as in case 1 of the proof of Theorem 6.4:
\itemitem{3.1} Compute a binomial
$\,b \in J \setminus P$.
\itemitem{3.2} Select an integer $d$
which might be sufficiently divisible.
\itemitem{3.3} Let $q$ be the largest
power of $char(k)$ that divides $d$ and set
$g := b^{[d]}/b^{[q]}$.
\itemitem{3.4} Compute a reduced
Gr\"obner basis ${\cal G}$ for the ideal  $(I:g)$.
\itemitem{3.5} If ${\cal G}$ consists of binomials,
call Algorithm 8.5 recursively with input ${\cal G}$.
\hfill \break
Otherwise return to step 3.2 and try a multiple of $d$.

\item{4.} If $L$ has infinite index in $L_\rho$ then
proceed as in case 2 of the proof of Theorem 6.4:
\itemitem{4.1} Compute a vector $m \in L_\rho$ whose image
in the quotient lattice $L\rho/L$ has infinite order.
Set $b := x^{m_+} - \rho(\sigma) x^{m_-}$.
\itemitem{4.2} Select an integer $d$
which might be sufficiently divisible.
\itemitem{4.3} Compute a reduced
Gr\"obner basis ${\cal G}$ for the ideal  $(I: b^{[d]})$.
\itemitem{4.5} If ${\cal G}$ consists of binomials,
call Algorithm 8.5 recursively with input ${\cal G}$.
\hfill \break
Otherwise return to step 4.2 and try a multiple of $d$.

\vskip .2cm

\noindent {\sl Comments: } The correctness of
Algorithm 8.6 follows from Theorem 6.4.

\vskip .4cm

\noindent {\bf Algorithm 8.7: Primary decomposition. }

 \noindent {\sl Input: \ } A cellular binomial ideal $I$.

\noindent {\sl Output: \ } Primary binomial ideals
$\, Q_i \, $ whose intersection is irredundant
and equals~$I$.

\item{1.} Compute the associated primes
$P_1,\ldots,P_s$ using Algorithm 8.5.

\item{2.} Choose a sufficiently large integer $e$.

\item{3.} For each prime $P_i$ do
\itemitem{3.1} If $char(k) = p > 0 $ then
let $\, R_i \,:=  \,I + P_i^{[p^e]} $.
\itemitem{3.2} If $char(k) = 0 $ then
let $\, R_i \,:=  \,I \,+\, M^e \,
+ \, (P_i \cap k[{\cal Z}]) $.
\itemitem{3.3} Compute  $Hull(R_i)$ using
Algorithm 8.6. Output $Q_i = Hull(R_i)$

\vskip .2cm

\noindent {\sl Comments: } The correctness of this
algorithm follows from Theorem 7.1'.
When computing with concrete binomial ideals,
it makes sense to replace $M^e$ in step 3.2
by $(x_i^{e_i},\, i \not\in {\cal Z})$ for sufficiently
large integers $e_i$. Good choices of these
integers, and many other algorithmic details,
will require further theoretical study and
experimentation.

\vskip .5cm

\noindent {\bf Examples 8.8.}
Here are a few examples of binomial primary
decompositions.
\item{(a)} The ideal
$\,I \,= \, (ab-cd, a^2, b^2, c^2, ac, bc) \,$
is primary but $I + (a)$ is not primary.
\item{(b)} The ideal $\, I \, = \,
( x^3 - y^3, x^4 y^5 - x^5 y^4) \,$ has the
following two primary decompositions:
$$ I \quad = \quad
(x-y) \,\cap\, \bigl( I + (x^9,y^9) \bigr)
\quad = \quad   (x-y) \cap
( x^2 + xy + y^2, \,
x^4 y^5 - x^5 y^4, x^{10},y^{10}).
$$
It can be shown that each primary decomposition of $I$
in which the embedded component has a quadratic
generator is {\sl not} binomial.
This proves that binomial ideals
behave differently from monomial ideals
with regard to the following result
of Bayer, Galligo and Stillman:
Every monomial ideal has a unique
``maximal primary decomposition''
in which each component is a monomial ideal
(see Eisenbud [1994], Exercise 3.17).
\item{(c)} The homogeneous ideal $\, I \,= \,
 ( c^5 - b^2 d^3, \, a^5 d^2 - b^7 ,\,
 b^5 - a^3 c^2 ,\, a^2 d^5 - c^7 ) \,$
is a circuit ideal (cf.~Example 7.9).
Its radical is the prime $\, P \, = \,I \,+ \,(ab-cd) $.
The projective toric variety defined by $P$
is the  rational normal curve of degree $7$.
The polyhedral cone ${\cal C}$ in formula (7.7)
has dimension $dim(P) = 2$. The faces of $P$
are $\{a,b,c,d\}$, $\,\{a\}$, $\,\{d\}$ and $\emptyset$.
The cellular decomposition has one component for each face:
$$ \eqalign{
I \quad = \quad P \,\,\,\,\cap \,\,\,\,
& (b^2 c^2 - a^2 d^2, \,
b^5-a^3 c^2 ,b^2 d^2, \,c^4 ,\, c^2 d^2 , \, d^4 )  \cr \,\,\cap
\,\,\, & ( b^2 c^2 - a^2 d^2 , \,c^5-b^2 d^3 ,\,
  a^2 c^2,\, b^4 ,\, a^2 b^2 ,\, a^4 )
\,\,\, \cap \,\,\, \bigl( I \, + \,
 (a^7,\, b^9,\, c^9,\, d^7 ) \bigr).\cr } $$
 This intersection is a primary decomposition of
$I$, as predicted by Theorem 7.6.
\item{(d)} The following radical binomial ideal appears
in (Eisenbud-Sturmfels [1993], Ex.~2.9),
$$ J \quad = \quad
( x_2 x_5 - x_1 x_6 , x_3, x_4) \,\,\cap \,\,
( x_1 x_4 - x_3 x_5 , x_2, x_6) \,\,\cap \,\,
( x_3 x_6 - x_2 x_4 , x_1, x_5) . $$
to show that the Noether complexity of an ideal
can be lower than that of any initial ideal.
Note that $J = I(\Delta)$ for a polyhedral complex $\Delta$
consisting of three quadrangles (cf.~Example 4.7 and Proposition 4.8).
It would be interesting to study the Noether complexity
of binomial ideals in general.
\item{(e)} We consider the typical
(but otherwise featureless) cellular binomial ideal
$$ \eqalign{
 I \quad  = \quad
\bigl( \,& bd^2-af^2, \,
    bce-acf,\,
    bcd-ace,\,
    b^2e-abf,\,
    b^2c,\,
    ae^2-bf^2,  \cr &
    ad^2-be^2,
    acd-bcf,
    abe-a^2f,
    abc,
    ab^2-b^3,
    a^2e-b^2f, \,
    a^2c,\, b^4,       \cr &
    a^2b-b^3,
    a^3-b^3,
    c^3e-c^3f,
    c^4,
    b^3d-b^3f,
    ac^3-bc^3,
    cd^4-ce^2f^2 \,\bigr)\cr} $$
It is cellular since
$ \sqrt{I}\, = \,(a,b,c )\,$ and
$d,e,f $ are non-zerodivisors mod~$I$.
 Using Algorithm 8.7 we may compute the
following primary decomposition:
$$ \eqalign{ & I \quad = \cr &
  (a,b,c) \quad \cap \quad \cr &
  (a,b,c^3,d^2-ef)\quad \cap \quad
  (a,b,c^3,d^2+ef)\quad \cap \quad  \cr &
  (a,b,c^4,e-f,d - i f)   \quad \cap \quad
  (a,b,c^4,e-f,d + i f)  \quad \cap \quad \cr &
  (a-b,b^4,c^4,b^2c,d-f,e-f)  \quad \cap \quad
  (a-b,b^2,c^3,bc,d+f,e+f)  \quad \cap \quad \cr &
  (a-b,b^3,c^4,bc,d+f,e-f)  \quad \cap \quad
  (a-b,b^2,c^3,bc,d-f,e+f)  \quad \cap \quad \cr &
  (a \! - \! \xi^2b,b^2,c^3,bc,d+\xi f, e+\xi^2f)   \,\, \cap \,\,
  (a \! - \! \xi^2b,b^3,c^3,bc,d \! - \! \xi f,e-\xi^2f)  \quad \cap  \cr &
  (a \! - \! \xi^2b,b^3,c^3,b^2c,d+\xi f,e\!- \!\xi^2f)  \,\,\, \cap \,\,\,
  (a\!- \! \xi^2b,b^2,c^3,bc,d\!-\!\xi f,e+\xi^2 f)  \,\,\, \cap  \cr &
  (a+\xi b,b^2,c^3,bc,d+\xi^2f,e-\xi f)   \,\,\, \cap \,\,\,
  (a+\xi b,b^3,c^3,b^2c,d \!-\! \xi^2f,e+\xi f)  \,\,\,   \cap \cr &
  (a+\xi b,b^3,c^3,bc,d+\xi^2f,e+\xi f)  \quad \cap \quad
  (a+\xi b,b^2,c^3,bc,d-\xi^2f,e-\xi f). \cr} $$
Here $i$ and $\xi$ are primitive roots of unity defined by
$\,i^2+1 =   \xi^2-\xi+1 = 0 $.
\item{(f)} The following binomial ideal appears in
(Koll\'ar, [1988], Example 2.3):
$$ (    x_1^{d_1},\,
    x_1 x_n^{d_2-1} \! - x_2^{d_2},\,
    x_2 x_n^{d_2-1} \! - x_3^{d_3},\,
\ldots, \,
x_{n-2} x_n^{d_{n-1}-1} \! - x_{n-1}^{d_{n-1}},
x_{n-1} x_n^{d_{n}-1} \! - x_{0}^{d_n} ).$$
This ideal has radical  $ (x_0,x_1,\ldots,x_{n-1})$ and it is primary.
This ideal provides a  lower bound for the effective Nullstellensatz
because it contains $x_0^{d_1 \cdots d_n}$
but not $x_0^{d_1 \cdots d_n-1}$.

\vskip 1cm

{
\baselineskip=12pt

\centerline {\bf References}

\vskip .1cm \par \frenchspacing

\item{} V.~I.~Arnold: A-graded algebras and continued fractions,
{\sl Communications in Pure and Appl. Math.} {\bf 42} (1989) 993-1000.
\vskip .1cm

\item{} D.~Bayer,  M.~Stillman:
On the complexity of computing syzygies.
{\sl Journal of Symbolic Computation} {\bf 6} (1988) 135-147.
\vskip .1cm

\item{} W.~D.~Brownawell:
Bounds for the degrees in the Nullstellensatz.
{\sl Annals of Math.} {\bf 126} (1987) 577-591.
\vskip .1cm

\item{} B.~Buchberger:
Gr\"obner bases - an algorithmic method in polynomial ideal theory,
Chapter 6 in N.K. Bose (ed.): {\sl Multidimensional Systems Theory},
D.~Reidel, 1985.
\vskip .1cm

\item{} D.~Buchsbaum and D.~Eisenbud:
Generic Free Resolutions and a Family of Generically Perfect Ideals.
{\sl Adv. in
Math.} {\bf 18} (1975) 245--301.
\vskip .1cm

\item{} D.~Buchsbaum and D.~S.~Rim: A generalized Koszul complex II:
depth and multiplicity.
{\sl Trans. Am. Math. Soc.} {\bf 111} (1964) 197--225.
\vskip .1cm

\item{} P.~Conti, C.~Traverso: Buchberger Algorithm
and Integer Programming, Proceedings AAECC-9 (New Orleans),
Springer Lect. Notes in Comp. Sci. {\bf 539} (1991) 130-139.
 \vskip .1cm

\item{} D.~Cox, J.~Little, D.~O'Shea:
{\sl Ideals, Varieties and Algorithms},
Springer, New York, 1992.
\vskip .1cm

\item{} W.~Decker, N.~Manolache, and F.~-O.~Schreyer:
Geometry of the Horrocks bundle on ${\bf P}^5$.
in {\sl Complex projective geometry (Trieste, 1989/Bergen,
1989)}, 128--148, London Math. Soc. Lecture Notes
179, Cambridge Univ. Press, Cambridge, 1992.
\vskip .1cm

\item{} P.~Diaconis and B.~Sturmfels:
Algebraic algorithms for generating from
conditional distributions, to appear.
\vskip .1cm

\item{} D.~Eisenbud: {\sl Commutative algebra with a view toward
algebraic geometry}. Springer-Verlag 1994 (to appear).
\vskip .1cm

\item{} D.~Eisenbud, C.~Hunecke, W.~Vasconcelos:
Direct methods for  primary decomposition. {\sl Inventiones
Math.} {\bf 110} (1992) 207--236.
\vskip .1cm

\item{} D.~Eisenbud, B.~Sturmfels: Finding sparse systems
of parameters. {\sl Journal of Pure and Applied Algebra}, to appear.
\vskip .1cm

\item{} W.~Fulton: {\sl Introduction to toric varieties}.
Annals of Math Studies 131, Princeton Univ. Press,
Princeton NJ (1993).
\vskip .1cm

\item{ } P.~Gianni, B.~Trager,  G.~Zacharias:
Gr\"obner bases and primary decomposition of
polynomial ideals, {\sl Journal of Symbolic Computation}
{\bf 6} (1988) 149--167.
\vskip .1cm

\item{} R.~Gilmer: {\sl Commutative semigroup rings}.
Chicago Lect. in Math.,
The University of Chicago Press, Chicago, 1984.
\vskip .1cm


\item{} J.~Koll\'ar: Sharp effective Nullstellensatz.
{\sl Journal Amer.~Math.~Soc.} {\bf 1},  (1988) 963--975.
\vskip .1cm


\item{} E.~Korkina, G.~Post, M.~Roelofs: Alg\`ebres gradu\'ees de type A.
{\sl Comptes Rendues de l'Acad.~de Sci.~Paris}, {\bf t.~314} ser. I num. 9,
(1992) 653-655.
\vskip .1cm

\item{} E.~Mayr, A.~Meyer:
The complexity of the word problem for commutative
semigroups and polynomial ideals.
{\sl Advances in Mathematics} {\bf 46} (1982) 305--329.
\vskip .1cm

\item{} I.~Hoveijn: Aspects of Resonance in Dynamical Systems,
Ph.D.~thesis, University of Utrecht, Netherlands, 1992.
\vskip .1cm

\item{} L.~Robbiano, M.~Sweedler: Subalgebra bases,
in W.~Bruns, A.~Simis (eds.): {\sl Commutative Algebra},
Springer Lecture Notes in Mathematics {\bf 1430}, 1990, pp.~61--87.
\vskip .1cm

\item{} R.~Stanley: Generalized H-vectors,
intersection cohomology of toric varieties, and related results.
in {\sl Commutative Algebra and Combinatorics,
Advanced Studies in Pure Math.} {\bf 11} (1987) 187--213.
\vskip .1cm

\item{} B.~Sturmfels:
Gr\"obner bases of toric varieties,
{\sl T\^ohoku Math. J.} {\bf 43} (1991) 249--261.
\vskip .1cm

\item{} B.~Sturmfels:
Asymptotic analysis of toric ideals,
{\sl Memoirs of the Faculty of
Sciences, Kyushu University, Series A: Mathematics}
{\bf 46}, No.~2, (1992) 217--228.
\vskip .1cm

\item{} R.~Thomas: A geometric Buchberger algorithm for
integer programming, to appear.
\vskip .1cm

\item{} S.~Xambo-Descamps: On projective varieties of minimal degree.
Collect. Math. {\bf 32} (1981) 149--163.
\vskip .1cm

\item{} C.K.~Yap: A new lower bound construction for
commutative Thue systems with applications,
{\sl Journal of Symbolic Computation} {\bf 12} (1991) 1--27.
\vskip .1cm

}

\bye